\documentclass{agtart_a}
\pdfoutput=1

\usepackage{graphicx,rotating}

\title{The universal Khovanov link homology theory}

\author{Gad Naot}
\givenname{Gad}
\surname{Naot}
\address{University of Toronto\\Department of Mathematics\\\newline
Toronto\\Ontario\\Canada}
\email{gnaot@math.toronto.edu}
\urladdr{}

\volumenumber{6}
\issuenumber{}
\publicationyear{2006}
\papernumber{65}
\startpage{1863}
\endpage{1892}

\doi{}
\MR{}
\Zbl{}


\received{19 May 2006}
\revised{14 July 2006}
\accepted{20 July 2006}
\published{1 November 2006}
\publishedonline{1 November 2006}
\proposed{}
\seconded{}
\corresponding{}
\editor{}
\version{}

\arxivreference{math.GT/0603347}

\let\xysavmatrix\xymatrix
\def\xymatrix{\disablesubscriptcorrection\xysavmatrix}
\AtBeginDocument{\let\bar\wbar\let\tilde\wtilde\let\hat\what}
\def\mod{\,\mathrm{mod}\,}
\makeatletter
\def\cnewtheorem#1[#2]#3{\newtheorem{#1}{#3}[section]
\expandafter\let\csname c@#1\endcsname\c@proposition}
\makeatother

\theoremstyle{plain}
\newtheorem{theorem}{Theorem}
\newtheorem{proposition}{Proposition}[section]
\cnewtheorem{lemma}[proposition]{Lemma}
\cnewtheorem{corollary}[proposition]{Corollary}

\theoremstyle{definition}
\cnewtheorem{definition}[proposition]{Definition}

\theoremstyle{remark}
\cnewtheorem{remark}[proposition]{Remark}

\newlength{\standardunitlength}
\newlength{\globalparindent}
\def\calF{{\mathcal F}}
\def\calO{{\mathcal O}}

\newcommand{\Cobl}{{\mathcal Cob}_{/l}}

\newcommand{\Komh}{\operatorname{Kom}_{/h}}
\newcommand{\Mat}{\operatorname{Mat}}
\newcommand{\Mor}{\operatorname{Mor}}

\newcommand{\eps}[2]{{\hspace{-3pt}\begin{array}{c}%
  \raisebox{-2.5pt}{\includegraphics[width=#1]{\figdir/#2}}%
\end{array}\hspace{-3pt}}}
\newcommand{\epsg}[2]{{\hspace{-3pt}\begin{array}{c}%
  \raisebox{0pt}{\includegraphics[width=#1]{\figdir/#2}}%
\end{array}\hspace{-3pt}}}

\begin{document}

\begin{asciiabstract}
We determine the algebraic structure underlying the geometric complex
associated to a link in Bar-Natan's geometric formalism of Khovanov's
link homology theory (n=2). We find an isomorphism of complexes
which reduces the complex to one in a simpler category.  This
reduction enables us to specify exactly the amount of information held
within the geometric complex and thus state precisely its universality
properties for link homology theories. We also determine its strength
as a link invariant relative to the different topological quantum
field theories (TQFTs) used to create link homology. We identify the
most general (universal) TQFT that can be used to create link homology
and find that it is ``smaller'' than the TQFT previously reported by
Khovanov as the universal link homology theory.  We give a new method
of extracting all other link homology theories (including Khovanov's
universal TQFT) directly from the universal geometric complex, along
with new homology theories that hold a controlled amount of
information. We achieve these goals by making a classification of
surfaces (with boundaries) modulo the 4TU/S/T relations, a process
involving the introduction of genus generating operators. These
operators enable us to explore the relation between the geometric
complex and its algebraic structure.
\end{asciiabstract}

\begin{htmlabstract}
We determine the algebraic structure underlying the geometric complex
associated to a link in Bar-Natan's geometric formalism of Khovanov's
link homology theory (n=2). We find an isomorphism of complexes
which reduces the complex to one in a simpler category.  This
reduction enables us to specify exactly the amount of information held
within the geometric complex and thus state precisely its universality
properties for link homology theories. We also determine its strength
as a link invariant relative to the different topological quantum
field theories (TQFTs) used to create link homology. We identify the
most general (universal) TQFT that can be used to create link homology
and find that it is &ldquo;smaller&rdquo; than the TQFT previously reported by
Khovanov as the universal link homology theory.  We give a new method
of extracting all other link homology theories (including Khovanov's
universal TQFT) directly from the universal geometric complex, along
with new homology theories that hold a controlled amount of
information. We achieve these goals by making a classification of
surfaces (with boundaries) modulo the 4TU/S/T relations, a process
involving the introduction of genus generating operators. These
operators enable us to explore the relation between the geometric
complex and its algebraic structure.
\end{htmlabstract}

\begin{abstract}
We determine the algebraic structure underlying the geometric complex
associated to a link in Bar-Natan's geometric formalism of Khovanov's
link homology theory ($n=2$). We find an isomorphism of complexes
which reduces the complex to one in a simpler category.  This
reduction enables us to specify exactly the amount of information held
within the geometric complex and thus state precisely its universality
properties for link homology theories. We also determine its strength
as a link invariant relative to the different topological quantum
field theories (TQFTs) used to create link homology. We identify the
most general (universal) TQFT that can be used to create link homology
and find that it is ``smaller'' than the TQFT previously reported by
Khovanov as the universal link homology theory.  We give a new method
of extracting all other link homology theories (including Khovanov's
universal TQFT) directly from the universal geometric complex, along
with new homology theories that hold a controlled amount of
information. We achieve these goals by making a classification of
surfaces (with boundaries) modulo the 4TU/S/T relations, a process
involving the introduction of genus generating operators. These
operators enable us to explore the relation between the geometric
complex and its algebraic structure.
\end{abstract}

\maketitle

\section{Introduction}

During the recent few years, starting with \cite{kho2}, Khovanov type
link homology theory has established itself as a dominant new field of
research within link invariant theory. Creating a homology theory
associated to each link, whose Euler characteristic is the Jones
polynomial, has proven itself to be a stronger invariant with many
advantages such as functorial properties regarding link
cobordisms. Together with the development of the algebraic language
used in the categorification process (the Khovanov type link homology
theories) there emerged a geometric/topological formalism describing
the entire process, due to Bar-Natan~\cite{ba1}. Initiating with the
task of clarifying \cite{kho2} and giving a visual geometric
description for ``standard'' Khovanov link homology, the geometric
formalism has evolved into a theory of its own. Using a fundamental
geometric language it was shown how to create an underlying framework
for Khovanov type link homology theories, which unifies, simplifies
and in many ways generalizes, much of the work done before.

\begin{sidewaysfigure}[p]
$\eps{180mm}{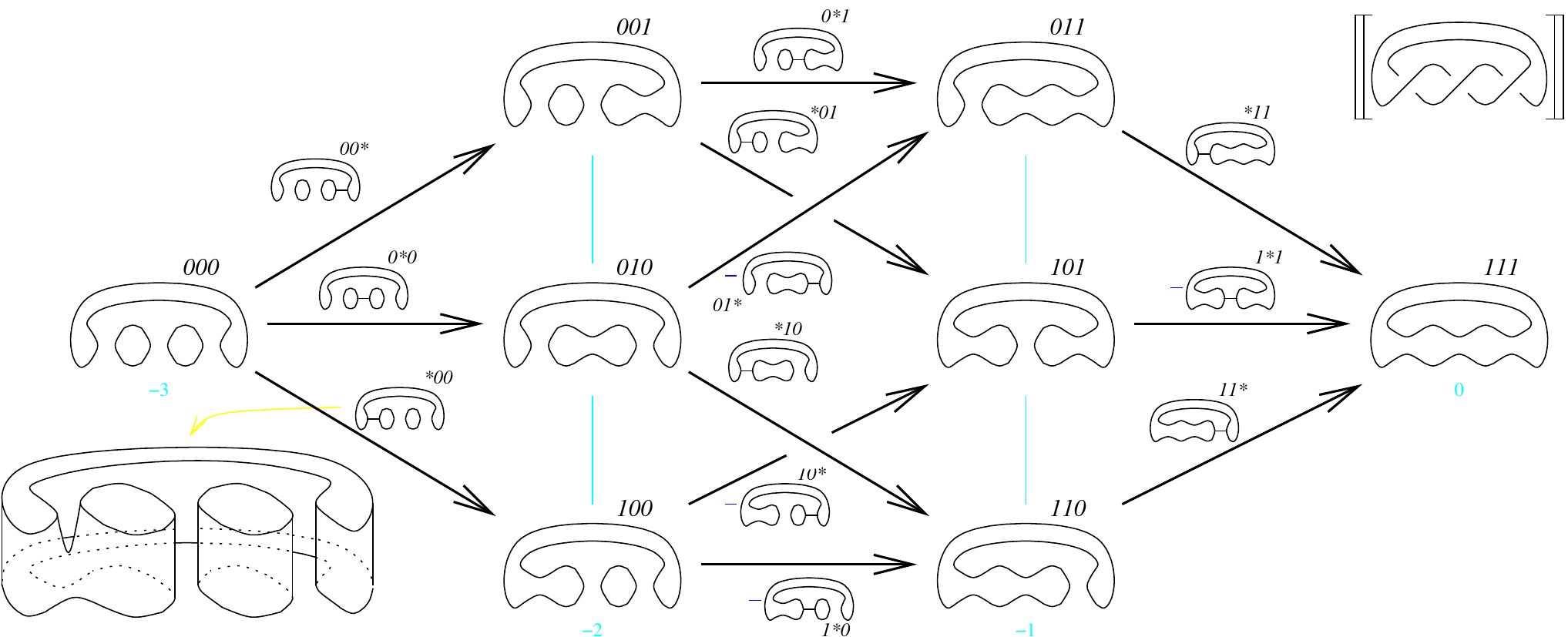}$
\caption{}\label{fig1}
\end{sidewaysfigure}

Though we assume familiarity with the main idea presented in
\cite{ba1}, the basic notion is as follows. Given a link diagram
$D$ one builds the ``cube of resolutions'' from it (a cube built of
all possible 0 and 1 smoothings of the crossings). The edges of the
cube are then given certain surfaces (cobordisms) attached to them
(with the appropriate signs). The entire cube is ``summed'' into a
complex (in the appropriate geometric category) while taking care of
some degree issues. We will call this complex the geometric complex
throughout the paper.  \fullref{fig1} should serve the
reader as a reminder of the process. The full description of it can
be found in \cite[Chapter 2]{ba1}.

The category in which one gets a link invariant is
$\Komh(\Mat(\Cobl))$, the category of complexes, up to homotopy,
built from columns and matrices of objects and morphisms
(respectively) taken from $\Cobl$, which is the category of
2--dimensional (orientable) cobordisms between 1 dimensional objects
(circles), where we allow formal sums of cobordisms over some ground
ring, modulo the following local relations:
\[ \text{The~4TU~relation:\qua}
\begin{array}{c}
\includegraphics[height=1cm]{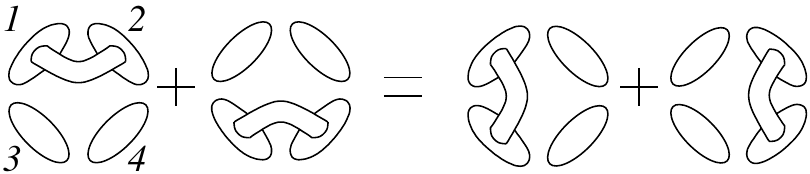}
\end{array}
\]
\[
\text{The~S~relation:\qua}\eps{7mm}{}=0
\]
\[
\text{The~T~relation:\qua}\eps{7mm}{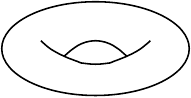}=2
\]
Using these relations a general theory was developed in which
invariance proofs became easier and more general, and homology
theories (TQFTs) became more natural --- one gets them by applying
``tautological'' functors on the geometric complex. Though studied
intensively in~\cite{ba1}, the full scope of the geometric theory
was not explored, and only various reduced cases (with extra
geometrical relations put and ground ring adjusted) were used in
connection with TQFTs and homology calculations. A full
understanding of the interplay between TQFTs used to create
different link homology theories and the underlying geometric
complex was not achieved although the research on the TQFT side is
considerably advanced, Khovanov~\cite{kho1}.

The objectives of this paper are to explore the full geometric
theory (working over $\mathbb{Z}$ with no extra relations imposed)
in order to answer the following questions:

\medskip
{\sl What is the \emph{algebraic} structure governing the full geometric
theory? How does the category $\Cobl$ look, and what can we do with
this information?

\medskip 
In what sense exactly is the geometric theory
universal? Does the full universal theory hold more information than
the different Khovanov type link homology theories (TQFTs) applied
to it?

\medskip 
What is the interplay between the different link
homology theories (TQFTs) and their geometric ``interpretation''
(via the geometric complex)?}

\medskip
We start by classifying surfaces with boundary modulo the 4TU/S/T
relations, and thus get hold of the way the underlying category of the
geometric complex looks. The main tools for this purpose will be using
\emph{genus generating operators} in order to extend the ground
ring. We prove a useful lemma regarding the free move of 2--handles
between components of surfaces in $\Cobl$ and using the genus
generating operators we get a simple classification of surfaces. We
present a reduction formula for surfaces in terms of a family of free
generators we identify.  The classification introduces the
topological/geometric motivation for the rest of the paper.
\fullref{sec2} presents the classification when 2 is invertible in the
ground ring and \fullref{sec4} presents the general case over $\mathbb{Z}$.

Then, we construct a reduction of the complex associated to a link.
We build an isomorphism of complexes and find that the complex
associated to a link is isomorphic to one that lives in a simpler
category. This category has only one object and the entire complex
is composed of columns of that single object. The complex maps are
matrices with monomial entries in one variable (a genus generating
operator). Thus the complex is equivalent to one built from free
modules over a polynomial ring in one variable. \fullref{sec3} presents
these results when 2 is invertible and \fullref{sec5} presents these
results for the general case over $\mathbb{Z}$.

The isomorphism of complexes gives us an immediate result regarding
the underlying algebraic structure of the universal geometric
theory. It presents us with a ``pre-TQFT'' structure of purely
topological/geometric nature. It turns out that the underlying
structure of the full geometric complex associated to a link (over
$\mathbb{Z}$) is the same as the one given by the \emph{co-reduced}
link homology theory using the following TQFT (over
$\mathbb{Z}[H]$):
\[
  \Delta_1: \begin{cases}
    v_+ \mapsto v_+\otimes v_- + v_-\otimes v_+  - H v_+\otimes v_+ &\\
    v_- \mapsto v_-\otimes v_- &
  \end{cases}
  \]
  \[
  m_1: \begin{cases}
    v_+\otimes v_-\mapsto v_- &
    v_+\otimes v_+\mapsto v_+ \\
    v_-\otimes v_+\mapsto v_- &
    v_-\otimes v_-\mapsto Hv_-
  \end{cases}
\]
In \fullref{sec6}, armed with the complex reduction above and the
underlying algebraic structure, we start exploring the most general
TQFTs that can be applied to the geometric complex to get a link
homology theory. This will result in a theorem which states that the
above co-reduced TQFT structure is the universal TQFT as far as
information in link homology is concerned.

We also get (\fullref{sec6}) a new procedure of extracting other TQFTs
directly from the geometric complex. This process is named
\emph{promotion} and can be used to get unfamiliar homology theories
which contain the information coming only from surfaces up to a
certain genus (in some sense a perturbation expansion in the genus),
thus extrapolating between the standard Khovanov link homology
theory and our universal one. This simplifies, completes and takes
into a new direction some of the results of~\cite{kho1}. Our paper
also generalizes some of the results of~\cite{ba1} and completes it
in some respects. We finish with some comments and further
discussion in \fullref{sec7}).

\subsubsection*{Acknowledgments}
This research was done as part of my PhD research at the University of
Toronto (Ontario, Canada).  I wish to thank Prof Bar-Natan for all his
support (and for allowing some figures to be borrowed).  I also wish
to thank the referee of this paper for his many useful suggestions and
comments.

\section{Classification of surfaces modulo the 4TU/S/T relations\\
when 2 is invertible}\label{sec2}

The geometric complex, an invariant of links and tangles, takes
values in the category $\Komh(\Mat(\Cobl))$. We wish to study this
category, as well as similar ones, in order to learn more about the
invariant. Reducing the geometric complex into a simpler complex in
a simpler category might give some further insight ---
computationally and theoretically. This section is devoted to the
study of the underlying category $\Cobl$ of 2--dimensional orientable
cobordisms between unions of circles in the simpler case when 2 is
invertible in the ground ring that we work over ($\mathbb{Q}$ or
$\mathbb{Z}[\frac{1}{2}]$, for example). Specifically we classify
all such surfaces modulo 4TU, S and T relations, which gives us the
morphism groups of $\Cobl$.

It is known~\cite{ba1} that when 2 is invertible in the ground ring
we work over, the 4TU relation is equivalent to \emph{the neck
cutting relation}:
\[
\text{NC ~relation:\qua}
2\begin{array}{c}\includegraphics[height=7mm]{}\end{array}
  =\begin{array}{c}\includegraphics[height=7mm]{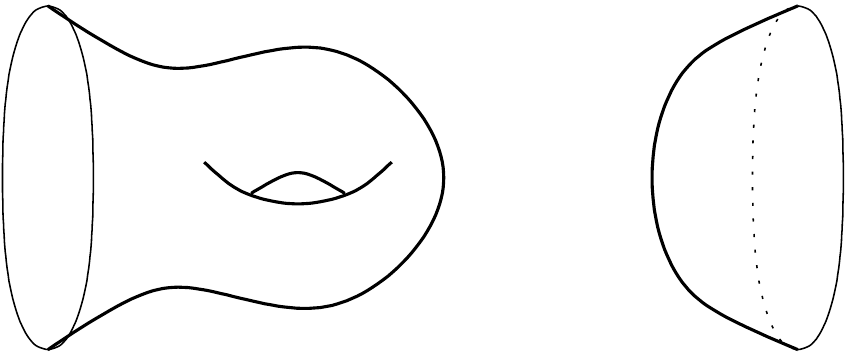}\end{array}
  +\begin{array}{c}\includegraphics[height=7mm]{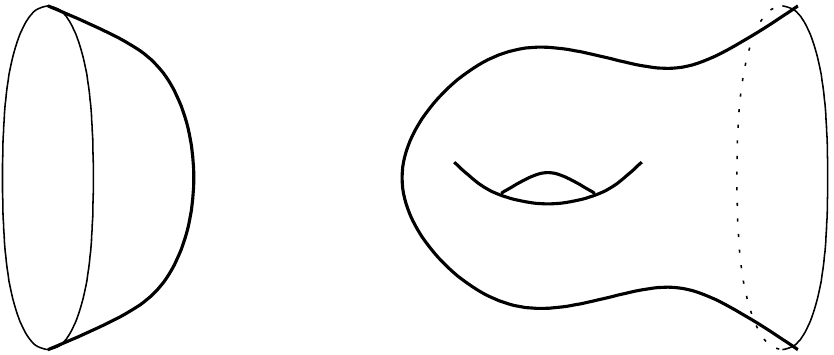}\end{array}
\]
We plan on using this simpler relation in this section.

\subsection{Notation}
For the sake of clarity and due to the nature of the topic, we will
try to give as many pictures and examples as possible. Still, one
needs some formal description from time to time, thus we will need
the use of some notation. Let $\Sigma_g(\alpha_1,\alpha_2,\cdots)$
denote a surface with genus g and boundary circles
$\alpha_1,\alpha_2,\cdots$. A disconnected union of such surfaces
will be denoted by
$\Sigma_{g_1}(\alpha_1,\cdots)\Sigma_{g_2}(\beta_1,\cdots)$. If the
genus or the boundary circles are not relevant for the argument at
hand, they will be omitted. Whenever we have a piece of surface
which looks like $\epsg{10mm}{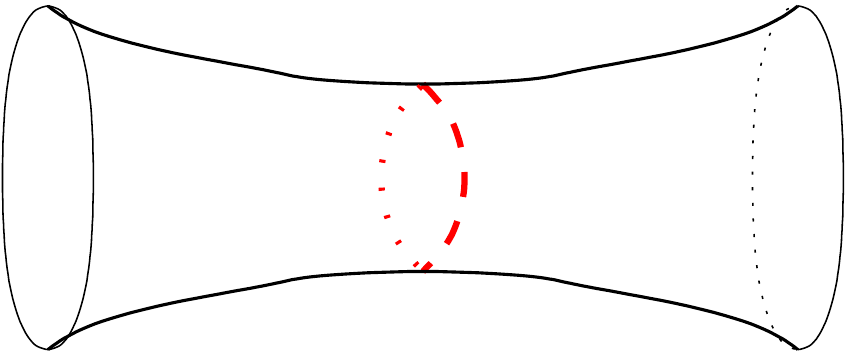}$ we will call it \emph{a neck}.
If cutting a neck separates the component into 2 disconnected
components then it will be a \emph{separating neck}, if not then it
is a \emph{non-separating neck} which means it is a part of \emph{a
handle} on the surface. A handle on the surface always looks locally
like $\eps{5mm}{}$, and by \emph{2--handle} on the surface we mean
a piece of surface which looks like $\eps{5mm}{}$. $\Cobl$ is
a general notation for either embedded 2 dimensional cobordism in 3
dimensional space (say a cylinder, like in~\cite{ba1}) or abstract
surfaces. One may choose whichever she likes, the theories (modulo
the relations) are the same (see \fullref{sec7}).

\subsection{The 2--handle lemma}
We start with proving a lemma that will become useful in classifying
surfaces modulo the 4TU/S/T relations.

\begin{lemma}
In $\Cobl$ 2--handles move freely between components of a surface.
In other words, modulo the 4TU relation, a surface with a 2--handle on one of
its connected components is equal to the same surface with the
2--handle removed and glued on a different component (see the picture
below the proof for an example).
\end{lemma}

\begin{proof}
The proof is an application of the neck cutting relation (NC), which
follows from the 4TU relation (over any ground ring). We look at a
piece of surface with a handle and a separating neck on it which
looks like $\eps{20mm}{}$. The rest of the surface continues
beneath the bottom circles and is not drawn. Applying the NC
relation to the vertical dashed neck and to the horizontal dashed
neck gives the following two equalities:
\[\eps{7mm}{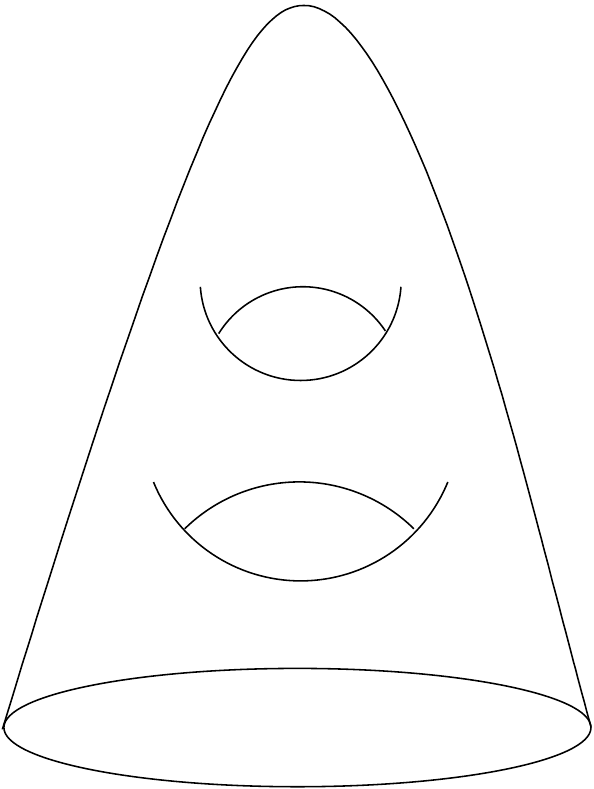}\eps{7mm}{} +
\eps{7mm}{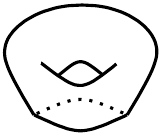}\eps{7mm}{vm} = 2\cdot\eps{15mm}{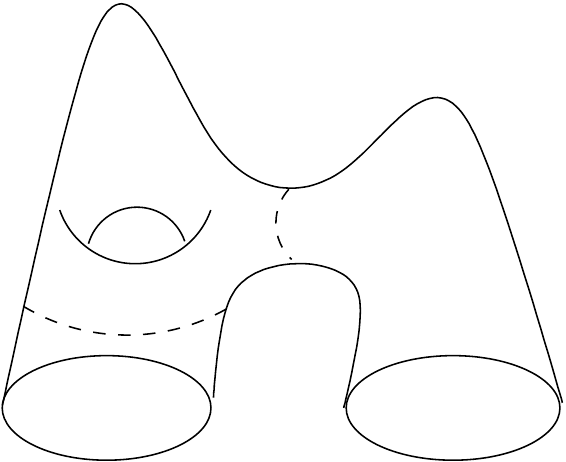} =
\eps{7mm}{vm}\eps{7mm}{vm} + \eps{7mm}{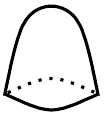}\eps{7mm}{2handle}
\]
We get $\eps{7mm}{2handle}\eps{7mm}{vp} =
\eps{7mm}{vp}\eps{7mm}{2handle}$. Since these are the top parts of
\emph{any} surface the lemma is proven.
\end{proof}

Notice that this lemma applies to any ground ring we
work over (in $\Cobl$), and does not require 2 to be invertible.

\medskip
{\bf Example}\qua The following equality holds in
$\Cobl$:
$$\eps{9mm}{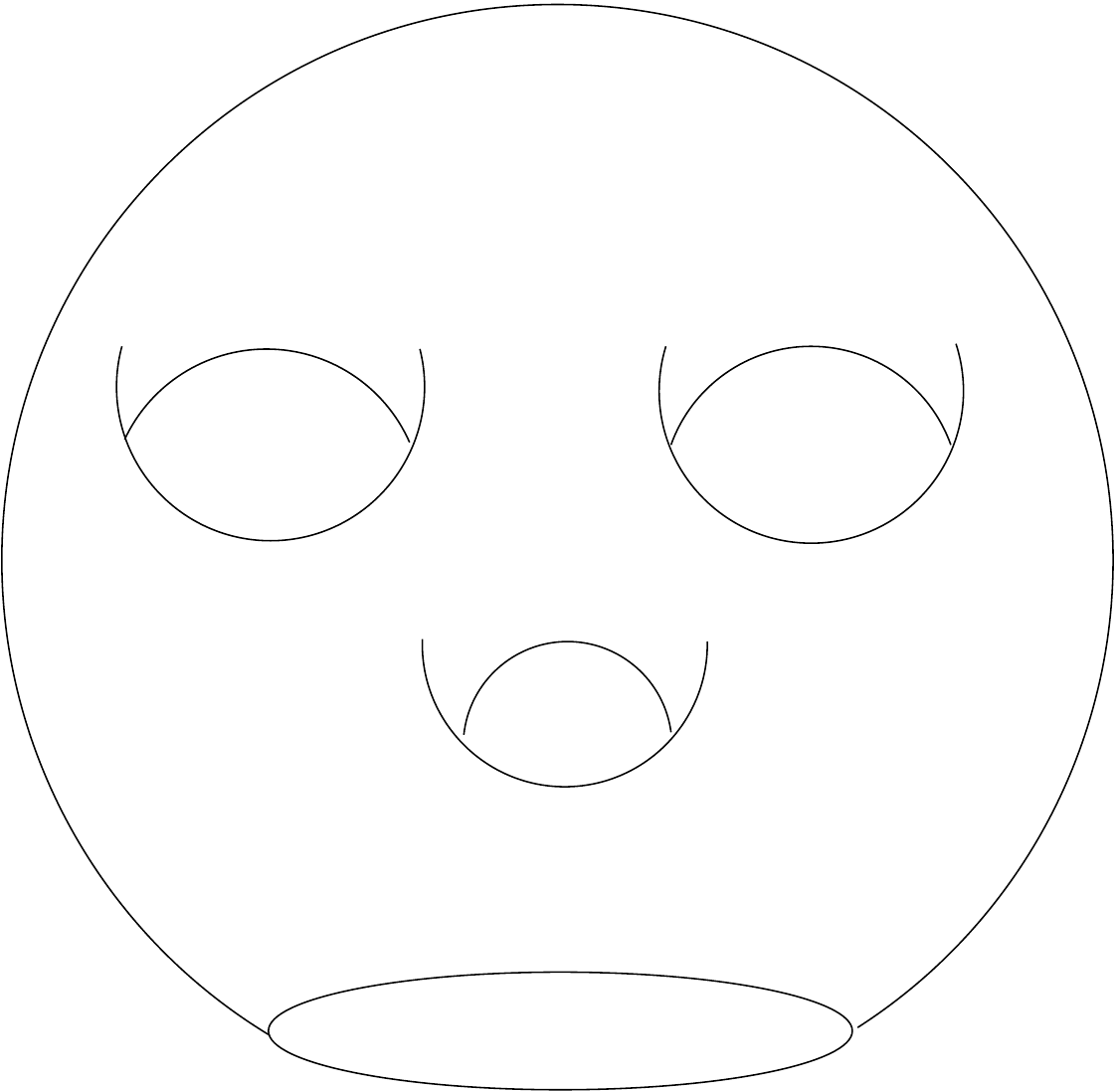}\eps{7mm}{2handle}=\eps{7mm}{vm}\eps{10mm}{}=
\eps{10mm}{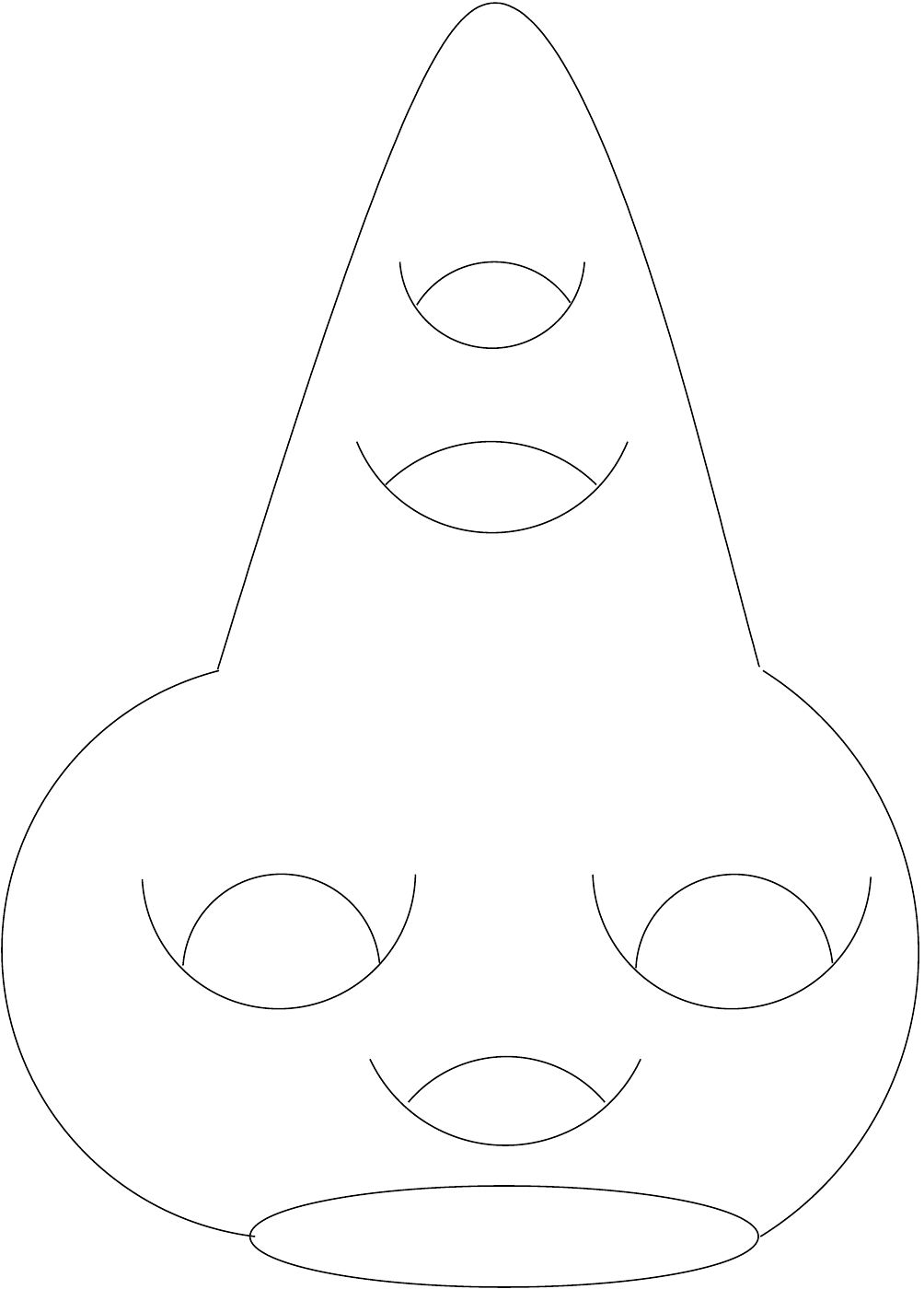}\eps{7mm}{vp}.$$

The 2--handle lemma allows us to reduce the classification problem to
a classification of surfaces with \emph{at most} one handle in each
connected component. This is done by extending the ground ring $R$
to $R[T]$ with $T$ a ``global'' operator (ie, acts \emph{anywhere}
on the surface) defined as follows:

\begin{definition}
\emph{The 2--handle operator}, denoted by $T$, is the operator that
glues a 2--handle somewhere on a surface (anywhere).
\end{definition}

\textbf{Examples}\qua $T\cdot\Sigma_{g_1}(\alpha)\Sigma_{g_2}(\beta) =
\Sigma_{g_1+2}(\alpha)\Sigma_{g_2}(\beta) =
\Sigma_{g_1}(\alpha)\Sigma_{g_2+2}(\beta)$ or $$T\cdot
\eps{6mm}{2handle}=\eps{10mm}{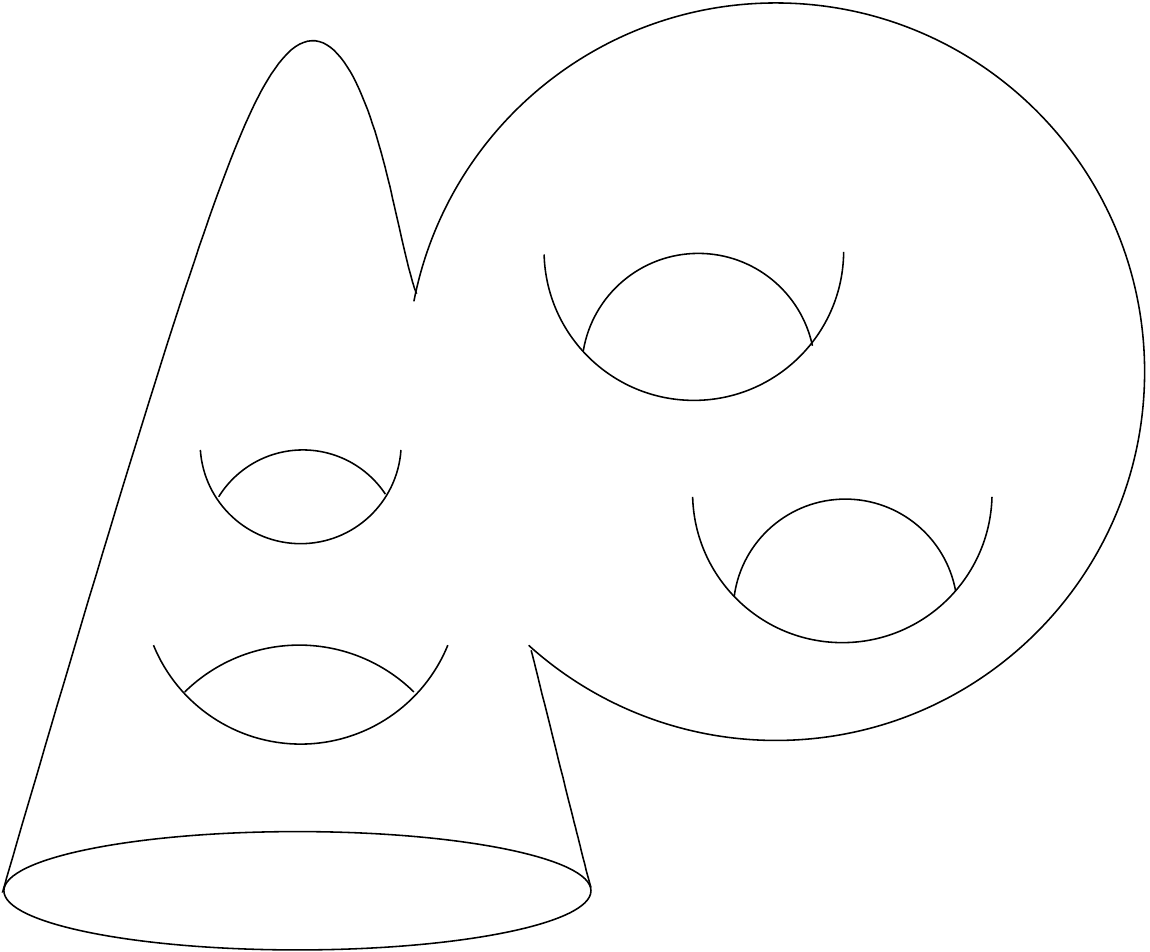}.$$

The 2--handle operator is given degree $-4$ according to~\cite{ba1}
degrees conventions, thus keeping all the elements of the theory
graded. Observe that since the torus equals 2 in $\Cobl$ (the T
relation), multiplying a surface $\Sigma$ with $2T$ is equal to
taking the union of $\Sigma$ with the genus 3 surface: \[2T \cdot
\Sigma = T \cdot \Sigma \cup \epsg{5mm}{} = \Sigma \cup
\epsg{9mm}{}\]
Note that the 2--handle operator does not operate on the empty
cobordism. Since 2 is invertible the operation of the 2--handle
operator can be naturally presented as the union with
$\epsg{9mm}{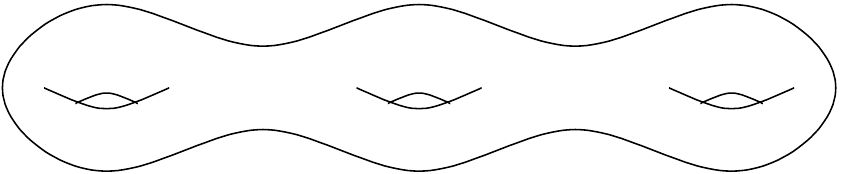}/2$, thus allowing meaningful operation on the empty
cobordism. We will adopt this presentation of $T$, and call it
\emph{the genus--3 presentation.}

\subsection{Classification of surfaces modulo 4TU/S/T over $\mathbb{Z}[\frac{1}{2}]$}

\begin{proposition}
Over the extended ring $\mathbb{Z}[\frac{1}{2},T]$, where $T$ is the
2--handle operator in the genus--3 presentation, the morphism groups
of $\Cobl$ (surfaces modulo S, T and 4TU relations) are generated
freely by unions of surfaces with exactly one boundary component and
genus 0 or 1 ($\eps{5mm}{vp},\eps{5mm}{vm}$), together with the
empty cobordism.
\end{proposition}

In other words, each morphism set from n circles to m circles is a
free module of rank $2^{n+m}$ over $\mathbb{Z}[\frac{1}{2},T]$,
where $T$ is the 2--handle operator in the genus--3 presentation.

\begin{proof}
Given a surface $\Sigma$ with (possibly) few connected components
and with (possibly) few boundary circles we can apply the 2--handle
lemma and reduce it to a surface that has at most genus 1 (i.e one
handle) in each connected component. This is done by extending the
ground ring we work over to $\mathbb{Z}[\frac{1}{2},T]$, where $T$
is the 2--handle operator.

When 2 is invertible the 4TU relation is equivalent to the NC
relation, and by dividing the relation equation by 2 we see that we
can cut any neck and replace it by the right hand side of the
equation:
$\begin{array}{c}\includegraphics[height=5mm]{\figdir/CNN}\end{array}
  =\frac{1}{2}\begin{array}{c}\includegraphics[height=5mm]{\figdir/CNL}\end{array}
  +\frac{1}{2}\begin{array}{c}\includegraphics[height=5mm]{\figdir/CNR}\end{array}$.
We will mod out this relation and reduce any surface into a
combination of free generating surfaces. We must make sure that the
reduction uses only the NC relation. Moreover the reduction must be
well defined (take the NC relation to zero), thus securing the
remaining generators from any relations.

Let us first look at the classification of surfaces with no boundary
at all. These are easy to classify due to the S and T relations.
Using the NC relation it is easy to show that $2\Sigma_{2g}=0$ and
with 2 invertible $\Sigma_{2g}=0$. Thus we can consistently extend
the ground ring with the operator $T$ and get: $\Sigma_g =
T^{\frac{g}{2}}\cdot\epsg{5mm}{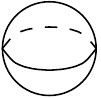}=0$ for $g$ even, and $\Sigma_g =
T^{\frac{g-1}{2}}\cdot\epsg{5mm}{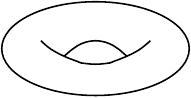}=2T^{\frac{g-1}{2}}$ for $g$
odd. Thus every closed surface is reduced to a ground ring element.
Allowing the empty cobordism (``empty surface''), every closed
surface will reduce to it over $\mathbb{Z}[\frac{1}{2},T]$ when we
use the genus--3 presentation of $T$. The closed surfaces morphism
set is thus isomorphic to $\mathbb{Z}[\frac{1}{2},T]$ and through
the genus--3 presentation acts on other morphism sets in the
category, creating a module structure.

Now, when cutting a separating neck in a component, one reduces a
component into two components at the expense of adding handles. Thus
using only the NC relation, over $\mathbb{Z}[\frac{1}{2},T]$, every
surface reduces to a sum of surfaces whose components contain
\emph{exactly} one boundary circle and at most one handle.

The reduction can be put in a formula form, considered as a map
$\Psi$ taking a surface $\Sigma_g(\alpha_1,\ldots,\alpha_n)$ into a
combination of generators (a sum of unions of surfaces with one
boundary circle and genus 0 or 1). This map is an isomorphism of
modules over the closed surfaces ring
($\mathbb{Z}[\frac{1}{2},T]$--modules):

\medskip
\textbf{The Surface Reduction Formula Over
$\mathbb{Z}[\frac{1}{2},T]$}
\begin{align}
\Psi \co  &\Sigma_g(\alpha_1,\ldots,\alpha_n)\\& \mapsto \frac{1}{2^{n-1}}
\cdot \sum_{i_1,\ldots,i_n=0,1}{(g+z)\mod2 \cdot
T^{\lfloor\frac{g+z}{2}\rfloor}
\cdot\Sigma_{i_1}(\alpha_1)\cdots\Sigma_{i_n}(\alpha_n)}\nonumber
\end{align}
where $z$ is the number of 0's among $\{i_1,\ldots,i_n\}$. If $n=0$
it is understood that the sum is dropped and replaced by $g\mod2
\cdot T^{\lfloor\frac{g}{2}\rfloor}$ times the empty surface. The
formula extends naturally to unions of surfaces.
\vspace{5pt}

A close look at the definition will show that indeed this map uses
only the NC relation to reduce the surface, in other words, the
difference between the reduced surface and the original surface is
only a combination of NC relations; (ie, $\Psi(\Sigma) - \Sigma =
\sum{NC}$). To see this, one need to show that there is at least one
way of using NC relations to get the formula. This can be done, for
instance, by taking the boundary circles one by one, and each time
cut the unique neck separating \emph{only} this circle from the
surface.
\vspace{5pt}

Last, we turn to the issue of possible relations between the
generators. These can be created whenever we apply $\Psi$ to two
surfaces which are related by NC relation and get two different
results. Thus the generators are free from relations if the
reduction formula is well defined, ie, respects the NC relation.
This is true, as the following argument proves $\Psi(NC)=0$. Look at
the NC relation, $\Sigma_{g_1+1}(\alpha)\Sigma_{g_2}(\beta) +
\Sigma_{g_1}(\alpha)\Sigma_{g_2+1}(\beta) -
2\Sigma_{g_1+g_2}(\alpha,\beta)$, where
$\alpha=\{\alpha_1,\ldots,\alpha_n\}$ and
$\beta=\{\beta_1,\ldots,\beta_m\}$ denote families of boundary
circles. Apply the formula everywhere:
\begin{align*}\frac{1}{2^{n+m-2}}&\cdot[(g_1+z_1+1\mod2)(g_2+z_2\mod2)
T^{\lfloor\frac{g_1+z_1+1}{2}\rfloor}
T^{\lfloor\frac{g_2+z_2}{2}\rfloor}\\& +(g_1+z_1\mod2)(g_2+z_2+1\mod2)
T^{\lfloor\frac{g_1+z_1}{2}\rfloor}T^{\lfloor\frac{g_2+z_2+1}{2}\rfloor}
\\&-
(g_1+g_2+z_1+z_2\mod2)T^{\lfloor\frac{g_1+z_1+g_2+z_2}{2}\rfloor}]\\&\cdot
\sum{\Sigma_{i_1}(\alpha_1)\cdots\Sigma_{i_n}(\alpha_n)\Sigma_{j_1}(\beta_1)\cdots\Sigma_{j_m}(\beta_m)}\end{align*}
where $z_1$ is the number of 0's among $\{i_1,\ldots,i_n\}$, $z_2$
is the number of 0's among $\{j_1,\ldots,j_m\}$ and the sum is over
these sets of indices as in formula (1). Going over all 4 options of
$(g_1+z_1,g_2+z_2)$ being (odd/even,odd/even) one can see the result
is always 0.
\vspace{5pt}

Formula (1) gives a map from surfaces to combinations of generators.
This map was shown to be well defined ($\Psi(NC)=0$), and uses only
the NC relations ($\Psi(\Sigma) - \Sigma = \sum{NC}$), thus
completing the proof.
\end{proof}

\section{Reduction of the geometric complex associated to a link and
the underlying algebraic structure over $\mathbb{Q}$}\label{sec3}

Armed with the classification of the morphisms of $\Cobl$ we are ready
to try and simplify the complex associated to a link in the case where
2 is invertible. We introduce an isomorphism of objects in $\Cobl$,
known as \emph{delooping}, that extends to an isomorphism in
$\Komh\Mat(\Cobl)$. This will reduce the complex associated to a link
into one taking values in a simpler category. One of the consequences
is uncovering the underlying algebraic structure of the geometric
complex. We will use the boxed surface notation where a boxed surface
has the geometric meaning of half a handle:
\[\eps{7mm}{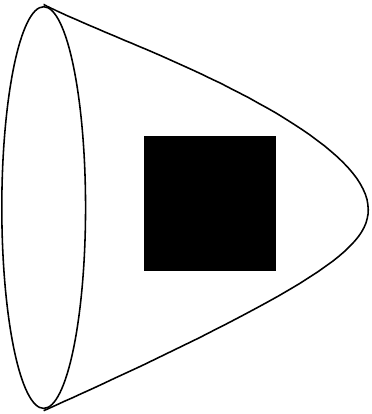}=\frac{1}{2}\eps{7mm}{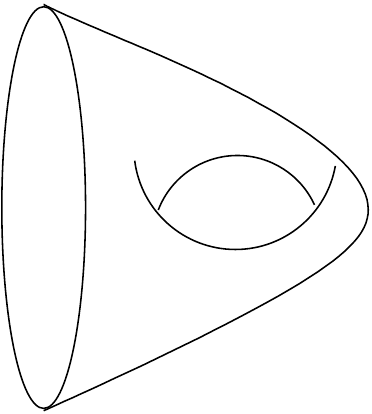}\]

\subsection{\textbf{The reduction theorem over $\mathbb{Q}$}}

\begin{theorem}\label{Thm1}
The complex associated to a link, over $\mathbb{Q}$ (or any ground
ring $R$ with 2 invertible), is equivalent to a complex built from
the category of free $\mathbb{Q}[T]$--modules ($R[T]$--modules
respectively).
\end{theorem}

This theorem is a corollary of the following proposition, which we
will prove first:

\begin{proposition}\label{Prop31}
In $\Mat(\Cobl)$, when 2 is invertible, the circle object is
isomorphic to a column of two empty set objects. This isomorphism
extends (taking degrees into account) to an isomorphism of complexes
in $\Komh\Mat(\Cobl)$, where the circle objects are replaced by
empty sets together with the appropriate induced complex maps.
\end{proposition}

\begin{definition}
The isomorphism that allows us to eliminate circles in exchange for
a column of empty sets is called \emph{delooping} and is given by:
\begin{center}
$\eps{20mm}{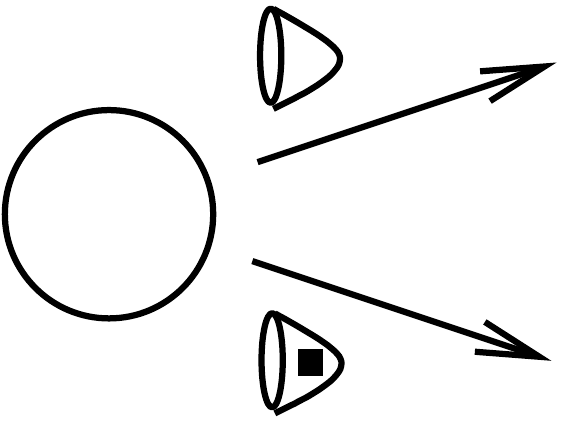}$ $ \left[ \begin{array}{c}
\emptyset \{-1\}\\
\\
\emptyset\{+1\}
\end{array} \right] $
$\eps{20mm}{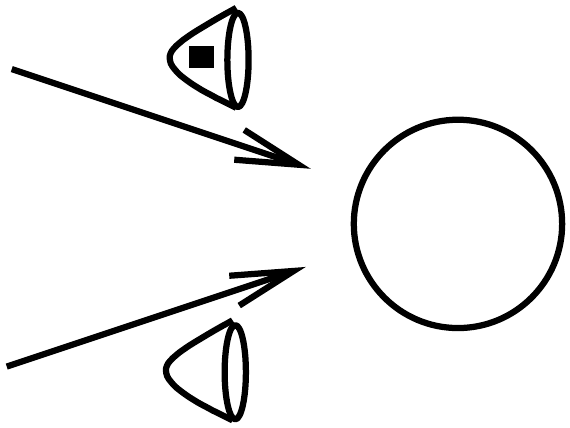}$
\end{center}
\end{definition}

\begin{proof}[Proof of \fullref{Prop31}]
It is easy to check that the compositions of the morphisms above are
the identity (thus making the circle isomorphic to the column of two
empty sets). One direction (from empty sets to empty sets) just
follows from the S and T relations (a boxed sphere equals 1), and
the other direction (from circle to circle) just follows from the NC
relation
$\begin{array}{c}\includegraphics[height=5mm]{\figdir/CNN}\end{array}
  =\frac{1}{2}\begin{array}{c}\includegraphics[height=5mm]{\figdir/CNL}\end{array}
  +\frac{1}{2}\begin{array}{c}\includegraphics[height=5mm]{\figdir/CNR}\end{array}$.

This proves that the circle object is isomorphic to a column of two
empty set objects in $\Mat(\Cobl)$. The extension to the category
$\Komh\Mat(\Cobl)$ is now straight forward by replacing every
appearance of a circle by a column of two empty sets (with the
degrees shifted as shown in the picture, maintaining the grading of
the theory). The induced maps are the composed maps of the
isomorphism above with the original maps of the complex.
\end{proof}

\begin{proof}[Proof of \fullref{Thm1}]
The theorem follows from the proposition. The only objects that are
left in the complex are the empty sets (in various degree shifts).
The induced morphisms in the reduced complex are just the set of
closed surfaces, isomorphic to $\mathbb{Q}[T]$ in the genus--3
presentation, thus making the complex equivalent to one in the
category of free $\mathbb{Q}[T]$--modules. Note that the isomorphism
in the proposition is between objects within the category
$\Mat(\Cobl)$ thus tautologically respects the 4TU/S/T relations.
One can check directly that the 4TU relation is respected --- since
the 4TU relation is a local relation in the interior of the surface
and the isomorphism is done on the boundary it obviously respects
it. One needs to check that the functor from $\Mat(\Cobl)$ to free
$\mathbb{Q}[T]$--modules is well defined (respects 4TU on the
morphism level) but this is obvious from the same reasons (it can be
checked directly using the classification from the previous
section). Thus the complex will look like columns of empty sets (the
single object) and maps which are matrices with entries in
$\mathbb{Q}[T]$. 
\end{proof}

{\bf Remark}\qua One can show that the morphisms
appearing in the complex associated to a link are further restricted
to matrices with $\mathbb{Z}[\frac{1}{2},T]$ entries. Furthermore,
we can work with boxed surfaces and re-normalize the 2--handle
operator by dividing it by 4, to get that the maps of the complex
associated to a link are always $\mathbb{Z}[T]$ matrices. Then,
degree considerations will tell us immediately that the entries are
only monomials in $T$.

\medskip
This theorem extends the simplification reported in
\cite{ba2,ba4}. There, it was done in order to compute the
``standard'' Khovanov homology, which in the context of the
geometric formalism means imposing one more relation which states
that any surface having genus higher than 1 is set equal to zero
(see also \cite[Section 9]{ba1}). In our context this means
setting the action of the operator $T$ to zero. When we set $T=0$
the complex reduces to columns of empty sets with maps being integer
matrices.

\subsection{\textbf{The algebraic structure underlying the complex
over $\mathbb{Q}$}} Now that we simplified the complex significantly
it is interesting to see what can we learn about the underlying
algebraic structure of the complex and about the complex maps in
terms of the empty sets.

As a first consequence we can see that the circle (the basic object
in the complex associated to a link) carries a structure composed of
two copies of another basic object (the empty set). It has to be
understood that this decomposition is a direct consequence of the
4TU/S/T relations and it is an intrinsic structure of the geometric
complex that reflects these relations. This already suggests that
any algebraic structure which respects 4TU/S/T and put on the circle
(for instance a Frobenius algebra, to give a TQFT) will factor into
a direct sum of two identical copies shifted by two degrees (two
copies of the ground ring of the algebra, say). This shows, perhaps
in a slightly more fundamental way, the result given in terms of
Frobenius extensions in~\cite{kho1}.

We denote $\emptyset\{-1\}$ by $v_-$ (comparing to~\cite{ba1}) or X
(comparing to~\cite{kho1}) and $\emptyset\{+1\}$ by $v_+$ or 1
(comparing accordingly). We also define a re-normalized 2--handle
operator $t$ to be equal the 2--handle operator $T$ divided by 4 (now
$t=\epsg{7mm}{g3}/8$ in the genus--3 presentation). Note that in the
category we are working right now this is just a notation, the
actual objects are the empty sets and we did not add any extra
algebraic structure. We use the tensor product symbol to denote
unions of empty sets and thus keeping track of degrees (remember
that all elements of the theory are graded).

Take the pair of pants map $\eps{7mm}{}$ between
$\bigcirc\bigcirc$ and $\bigcirc$, and look at the two complexes
$\bigcirc\bigcirc\xrightarrow{\rm pants}\bigcirc$ and
$\bigcirc\xrightarrow{\rm pants}\bigcirc\bigcirc$. Reduce these
complexes using our theorem (ie, replace the circles with empty
sets columns and replace the complex maps with the induced maps). We
denote by $\Delta_2$ the following composition:

\begin{center}
$ \Delta_2: \left[ \begin{array}{c}
v_-\\
v_+
\end{array} \right]
\xrightarrow{\rm isomorphism} \bigcirc
\xrightarrow{\eps{7mm}{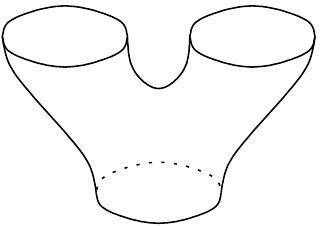}} \bigcirc\bigcirc
\xrightarrow{\rm isomorphism} \left[
\begin{array}{c}
v_-\otimes v_-\\
v_-\otimes v_+\\
v_+ \otimes v_-\\
v_+ \otimes v_+
\end{array} \right]
$
\end{center}

We denote by $m_2$ the composition in the reverse direction. By
composing surfaces and using the re-normalized 2--handle operator it
is not hard to see that the maps we get are:
\[
  \Delta_2: \begin{cases}
    v_+ \mapsto \left[ \begin{array}{c} v_+\otimes v_- \\ v_-\otimes v_+ \end{array} \right] & \\
    v_- \mapsto \left[ \begin{array}{c} v_-\otimes v_- \\ t v_+\otimes v_+ \end{array} \right] &
  \end{cases}
  \qquad
  m_2: \begin{cases}
    v_+\otimes v_-\mapsto v_- &
    v_+\otimes v_+\mapsto v_+ \\
    v_-\otimes v_+\mapsto v_- &
    v_-\otimes v_-\mapsto tv_+.
  \end{cases}
\]
The maps have no linear structure, these are just
$2\times4$ matrices of morphisms in the geometrical category.
Nonetheless this is exactly the pre-algebraic structure of the
generalized Lee TQFT. See~\cite{lee1} for the non generalized one
(where $t=1$) and \cite[section 9]{ba1}. The above TQFT is
denoted $\calF_3$ in~\cite{kho1}. It is important to note that so
far in our category we did not apply any TQFT or any other functor
to get this structure, these maps are not multiplication or
co-multiplication in an algebra. It is all done intrinsically within
our category, and as we will see later it imposes restrictions on
the most general type of functors one can apply to the geometric
complex. One can say that the geometric structure over $\mathbb{Q}$
has the underlying structure of this specific \emph{``pre-TQFT''}.

\section{Classification of surfaces modulo 4TU/S/T relations over $\mathbb{Z}$} 
\label{sec4}

When 2 is not invertible the neck cutting relation is not equivalent
to the 4TU relation and we need to make sure we use the ``full
version'' of the 4TU relation in reducing surfaces into generators.
Still there is a simpler equivalent version of the 4TU relation which
involves only 3 sites on the surface:
\begin{eqnarray*}
  3S_1:\qquad& \eps{7cm}{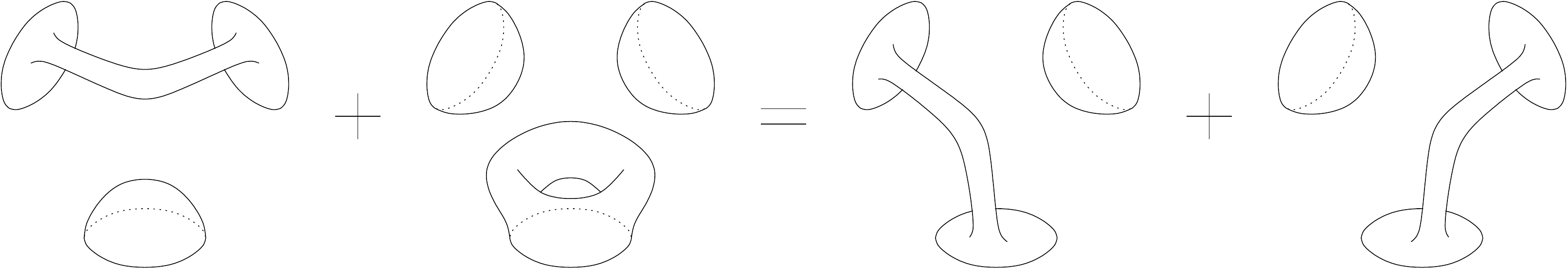} \\
  3S_2:\qquad& \displaystyle
  \sum_{\parbox{2cm}{\begin{center}\scriptsize\rm
    $0^\circ$, $120^\circ$, $240^\circ$ \newline
    rotations
  \end{center}}}
  \left(\eps{1.2cm}{}-\eps{1.2cm}{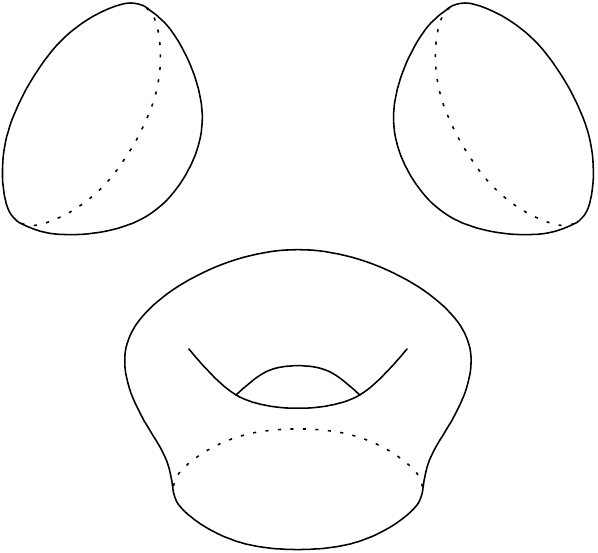}\right) = 0
\end{eqnarray*}
The version that will be most convenient for us is the $3S1$
relation. We will use the 2--handle lemma (which does not depend on
the invertibility of 2), extend the ground ring again, and identify
a complete set of generators.

\begin{definition}
Assume a surface has at least one boundary circle. Choose one
boundary circle. The component that contains the specially chosen
boundary circle is called \emph{the special component}. All other
components will be called \emph{non-special}.
\end{definition}

\begin{definition}
The \emph{special 1--handle operator}, denoted $H$, is the operator
that adds a handle to the special component of a surface.
\end{definition}

Note that $T=H^2$, though $H$ acts ``locally'' (acts \emph{only} on
the special component) and $T$ acts ``globally'' (anywhere on the
surface). We extend our ground ring to $\mathbb{Z}[H]$. We remind
the reader that the entire theory is still graded with $H$ given
degree -2.

\begin{remark}
From now on we will assume that all surface has at least one
boundary component, thus there is always a special component, and
the action of $H$ is well defined (by choosing a special component).
\end{remark}

\begin{proposition}
Over the extended ring $\mathbb{Z}[H]$ (after a choice of special
boundary circle) the morphisms of $\Cobl$ with one of the
source/target objects non empty (ie, surfaces modulo 4TU/S/T
relations with at least one boundary circle) are generated freely by
surfaces which are composed of a genus zero special component with
any number of boundary circles on it and zero genus non-special
components with exactly one boundary circle on them. $H$ is the
special 1--handle operator.
\end{proposition}

An example for a generator would be $\eps{20mm}{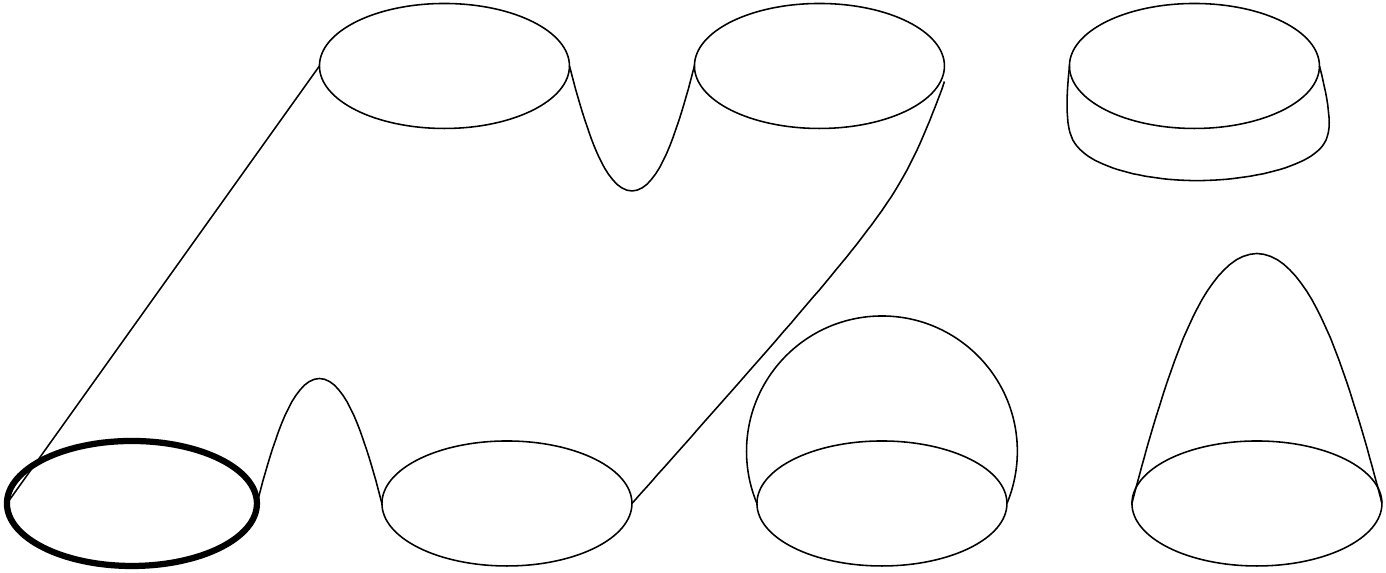}$ (the special
circle is at the bottom left).

Another example would be the ``Shrek surface'':
$\eps{20mm}{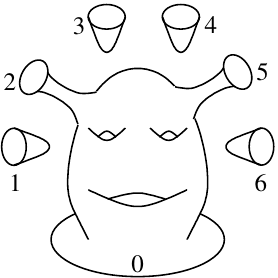}$. This surface has a special boundary circle
marked with the number zero (Shrek's neck), 6 other boundary circles
(number 2 and 5 belong to the special component --- the head) and 3
handles on the special component. Over the extended ring, this
surface is generated by the ``Shrek shadow'' $\eps{20mm}{}$
and equals to $H^3\eps{20mm}{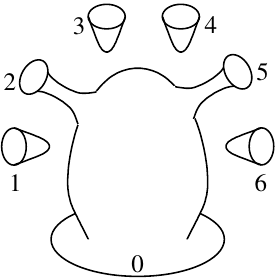}$.

\begin{remark}
In order for the composition of surfaces (morphisms) to be
$H$--linear, and thus respect the $\mathbb{Z}[H]$--module structure,
one needs to restrict to a category where all the morphisms preserve
the special component, ie, we always have the special boundary
circles of the two composed morphisms in the same component. As we
will see (next section) this happens naturally in the context of
link homology.
\end{remark}

The rest of this section is devoted to the proof of the proposition
above.

\begin{proof}
We start with a couple of toy models. If the surface components
involved in the $3S2$ (or $3S1$) relation has all together one
boundary circle, then the $3S2$ (or $3S1$) relation is equivalent to
NC relation, and further more up to the 2--handle lemma it is
trivially satisfied. This makes the classification of surfaces with
only one boundary component easy. Extend the ground ring to
$\mathbb{Z}[T]$ and the generators will be $\eps{0.5cm}{vm}$ and
$\eps{0.5cm}{vp}$. We can now go one step further and use the
``local'' 1--handle operator $H$. Over $\mathbb{Z}[H]$ there is only
one generator $\eps{0.5cm}{vp}$.

If the surface's components involved have only two boundary circles
all together, then the $3S2$ (or $3S1$) relation is again equivalent
to the NC relation. We can classify now all surfaces with 2 boundary
circles modulo the 4TU relation over $\mathbb{Z}$. Let us extend the
ground ring to $\mathbb{Z}[T]$ thus reducing to components with
genus 1 at most. Then, we choose one of the boundary circles (call
it $\alpha$), and use the NC relation to move a handle from the
component containing the other boundary circle (call it $\beta$) to
the component containing $\alpha$ (ie,
$\Sigma_1(\beta)\Sigma_g(\alpha) = 2\Sigma_g(\beta,\alpha) -
\Sigma_0(\beta)\Sigma_{g+1}(\alpha)$). This leaves us with the
following generators over $\mathbb{Z}[T]$: $\eps{1cm}{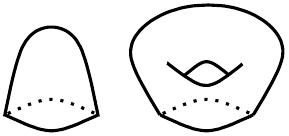}$,
$\eps{5mm}{vp}\eps{5mm}{vp}$, $\eps{10mm}{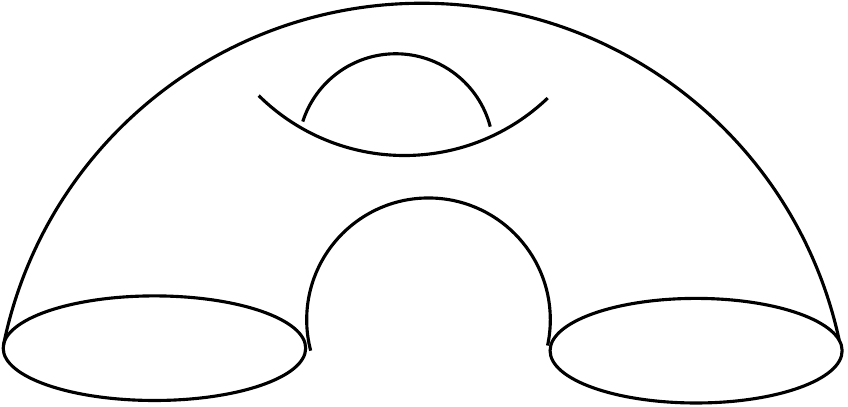}$ and
$\eps{1cm}{}$ ($\alpha$ is drawn on the right). This asymmetry
in the way the generating set looks like is caused by the symmetry
braking in the way the NC relation is applied (we chose the special
circle $\alpha$, on the right). On the other hand it allows us again
to replace the ``global'' 2--handle operator $T$, with the ``local''
1--handle operator $H$, which operates only on the component
containing $\alpha$. Thus, over $\mathbb{Z}[H]$ we only have the
following generators: $\eps{5mm}{vp}\eps{5mm}{vp}$ and
$\eps{1cm}{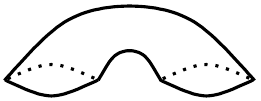}$. This set is symmetric again, and the asymmetry
hides within the definition of $H$.

The place where the difference between the NC and $3S2$ (or $3S1$)
relations comes into play is when three different components with
boundary on each are involved. We would like to use, in the general
case, the $3S1$ relation as a neck cutting relation for the neck in
the upper part of $\eps{10mm}{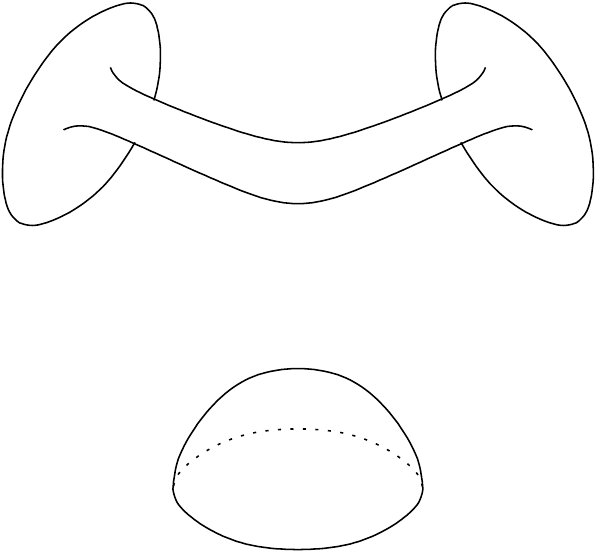}$ (by replacing it with the other
3 elements in the equation). The way the relation is applied is not
symmetric (not even visually), while the upper 2 sites are
interchangeable the lower site plays a special role. This hints that
this scheme of classification would be easier when choosing a
special component of the surface. This is done by choosing a special
boundary circle as in definition 4.1.

\textbf{Cutting non separating necks over $\mathbb{Z}$}\qua Assume
first that we only use $\eps{10mm}{3S2}$ in the $3S1$ relation to
cut a non-separating neck (the upper two sites remain in one
connected component after the cut), and we take the lower site to be
on the special component. Now, whenever we have a handle in a
non-special component (dashed line on the left component in the
picture below) we can move it using the $3S1$ relation to the
special component (the component on the right) at the price of
adding surfaces that has the non-special component connected to the
special component, thus becoming special: $\eps{20mm}{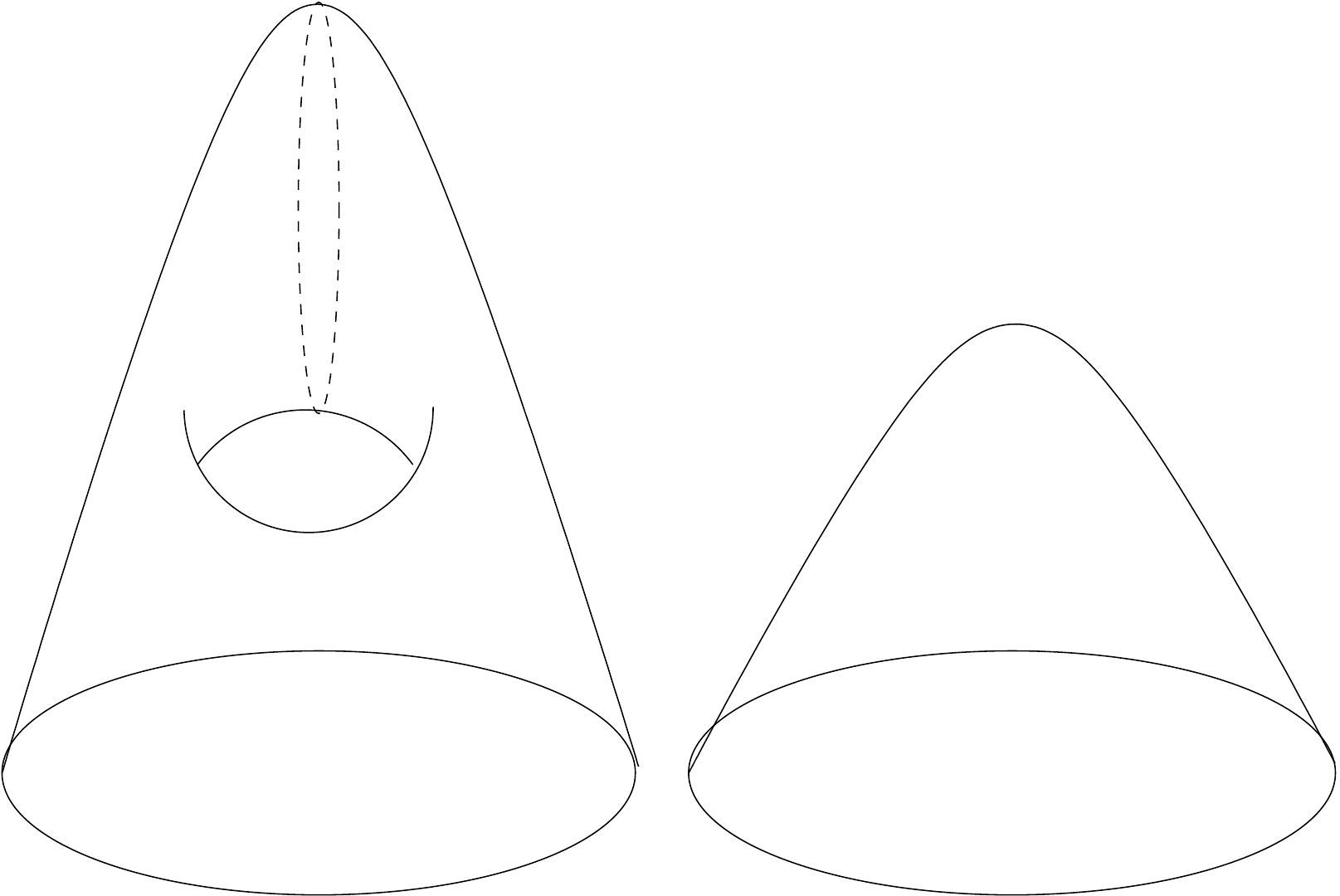} =
\eps{20mm}{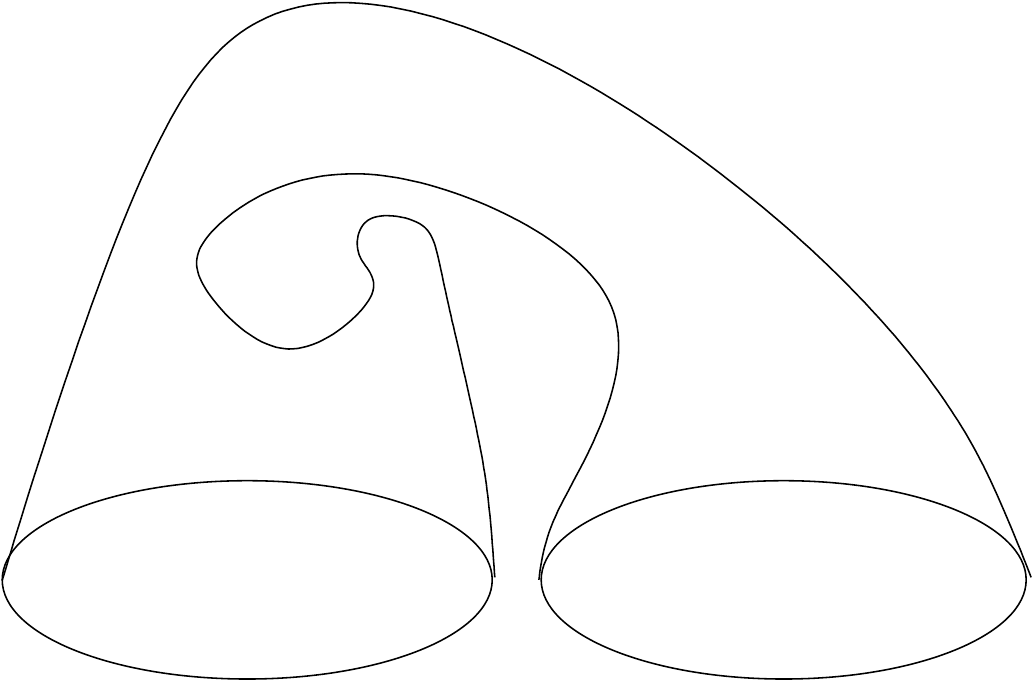} + \eps{20mm}{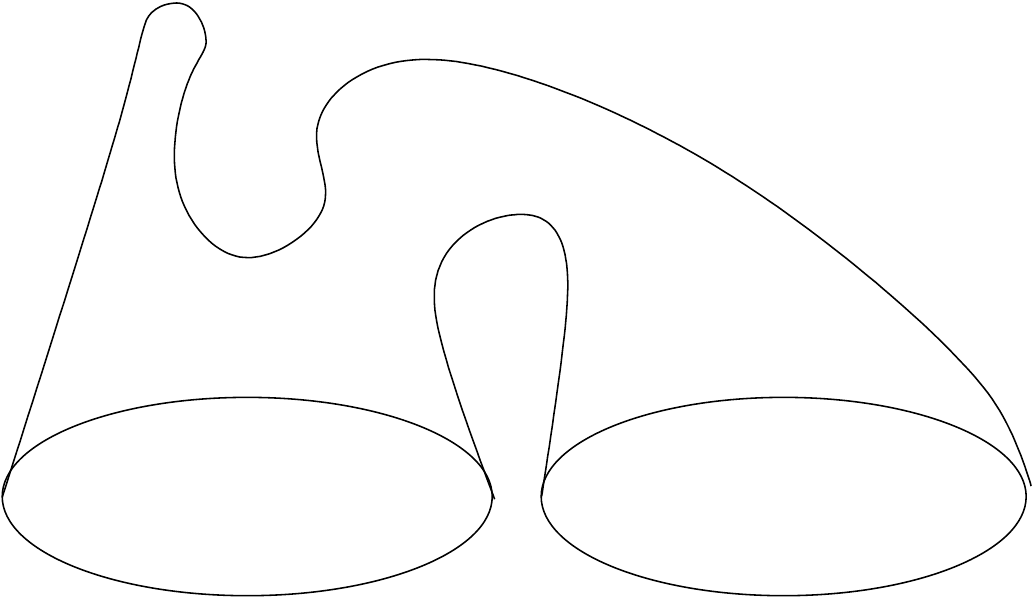} -
\eps{20mm}{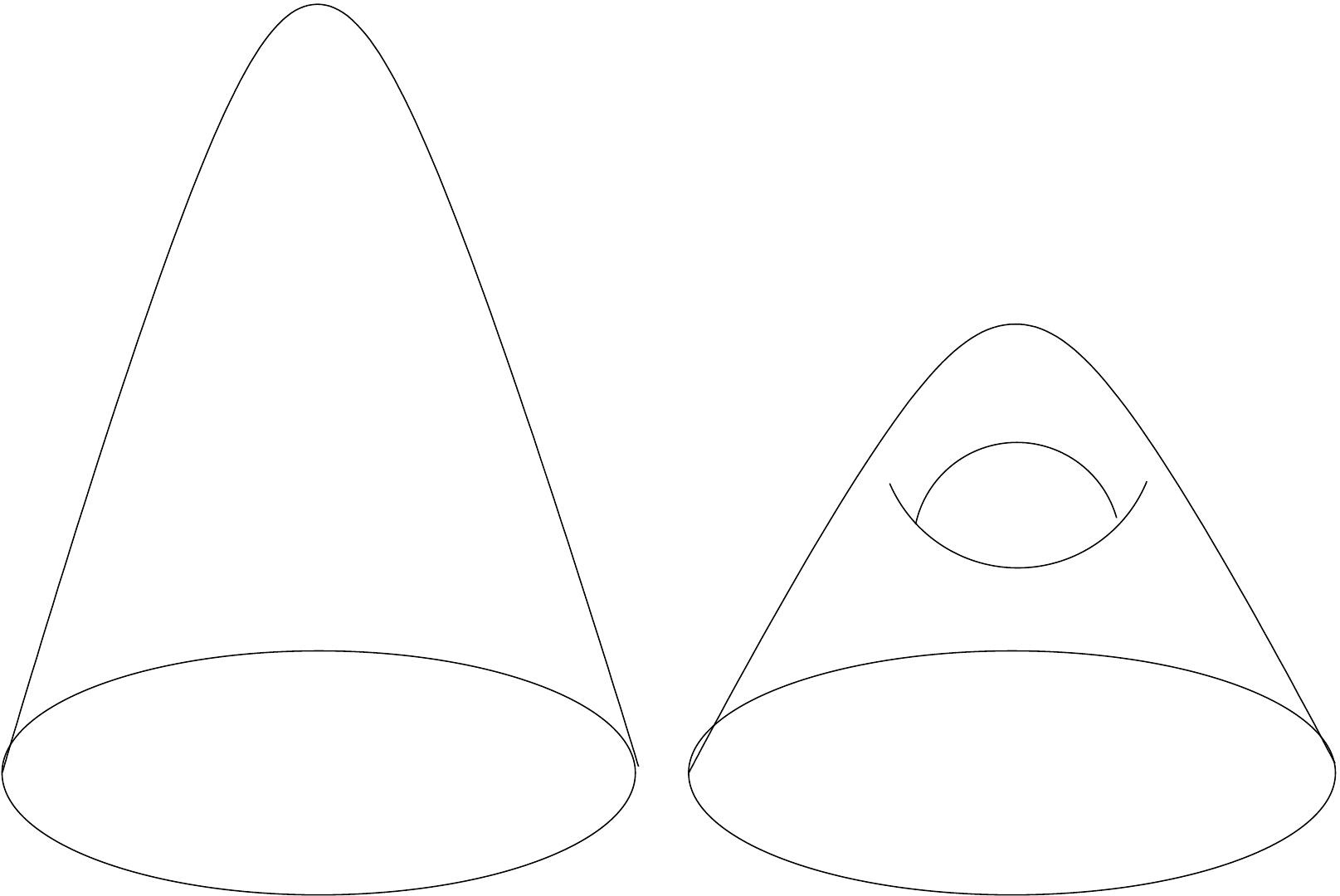}$. After this reduction we are left with
surfaces involving handles only on the special component, and the
rest of the components are of genus zero (with any amount of
boundary circles on them). For example $\eps{20mm}{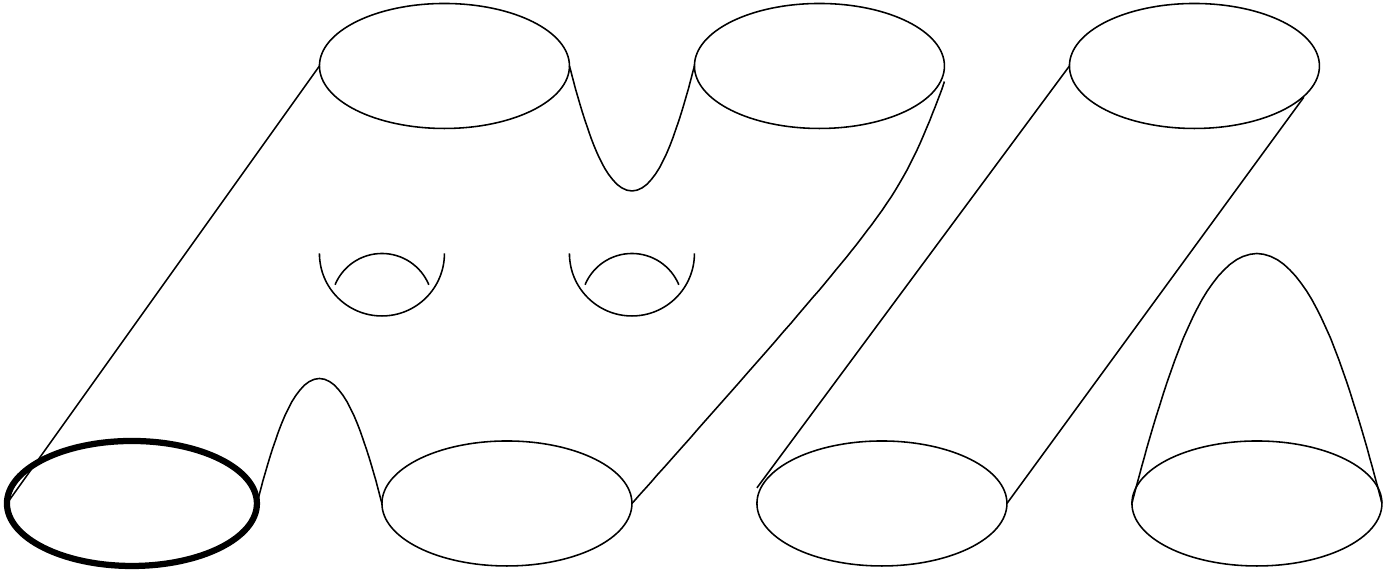}$, where the
special circle is at the bottom left.
\vspace{7pt}

\textbf{Cutting separating necks over $\mathbb{Z}$}\qua We reached a
point where the surfaces left are unions of a special component with
any number of handles on it and zero genus non-special components
with any possible number of boundary circles. We can still use the
$3S1$ relation to cut separating necks (the two upper sites in
$\eps{10mm}{3S2}$ will be on two disconnected components after the
cut). We will use the relation with the upper two sites on a
non-special component and the lower site on the special component.
The result is separating the non-special component into two at the
price of adding a handle to the special component and adding
surfaces that have connecting necks between the special component
and what used to be the non-special component (the other summands of
$3S1$). Thus, whenever we have a non-special component that has more
than one boundary circle, we can reduce it to a sum of components
that are either special or have exactly one boundary circle. An
example for a summand in a reduced surface is $\eps{20mm}{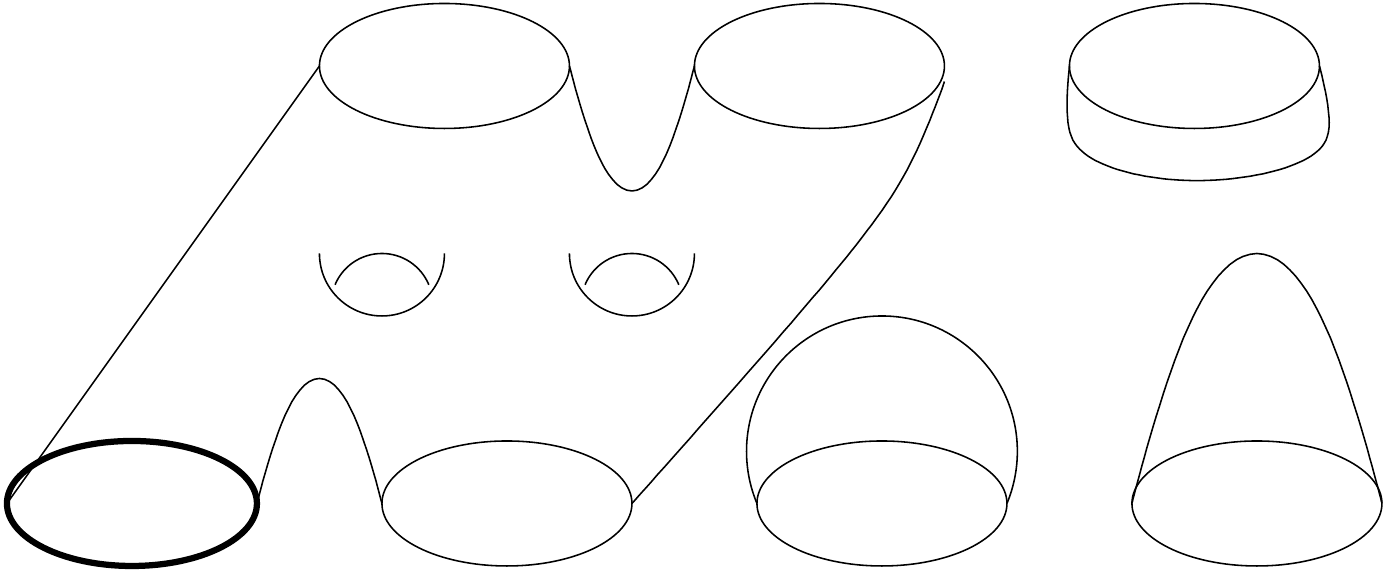}$,
where the special circle is at the bottom left.
\vspace{7pt}

\textbf{Ground ring extension, surface reduction formula over
$\mathbb{Z}$}\qua We extend our ground ring to $\mathbb{Z}[H]$. We will
put all the above discussion into one formula for reducing any
surface into generators (considered as a map $\Psi$). Let $\alpha$
denote a family of boundary circles $\{\alpha_1,\ldots,\alpha_n\}$
and let $S$ denote the special boundary circle. We will use
$\beta\in 2^\alpha$ for any subset
$\{\alpha_{i_1},\ldots,\alpha_{i_k}\}$ of $\alpha$, and
$\{\alpha_{i_{k+1}},\ldots,\alpha_{i_n}\}$ for its complement in
$\alpha$.

\medskip \textbf{The surface reduction formula over $\mathbb{Z}$}

{Special component reduction:}
\begin{equation}
\Sigma_g(S,\ldots)=H^g\cdot\Sigma_0(S,\ldots)
\end{equation}

{Non-special component genus reduction ($g\geq1$):}
\begin{equation}
\Sigma_0(S,\ldots)\Sigma_g(\alpha) = 2\cdot(g\mod2)\cdot H^{g-1}
\Sigma_0(S,\ldots,\alpha) + (-1)^g H^g
\Sigma_0(S,\ldots)\Sigma_0(\alpha)
\end{equation}

{Non-special component neck reduction ($n\geq2$):}
\begin{equation}
\Sigma_0(S,\ldots)\Sigma_0(\alpha) =
\sum_{\beta\in2^\alpha,\beta\neq\alpha}{(-1)^{n-k-1}H^{n-k-1}
\Sigma(S,\ldots,\beta) \Sigma_0(\alpha_{i_{k+1}}) \cdots
\Sigma_0(\alpha_{i_n})}
\end{equation}
Formula (4) is plugged into the result of formula (3)
in order to reduce any non-special component into combinations of
generators. Extend the above formulas to unions of non-special
components in a natural way (just iterate the use of the formula)
and get a reduction of any surface into generators. The previous
discussion shows that this formula uses the $3S1$ relation only
(ie, $\Psi(\Sigma)-\Sigma=\sum{3S1}$).

Finalizing the proof is done by taking the formula and checking that
it is a map from surfaces to the generators (by definition) that
uses only the 4TU ($3S1$) relations (proved by the discussion above)
and that the generators are free (ie, it is well defined). The last
part can be done by a direct computation (applying the formula to
all sides of a 4TU relation to get zero) or by applying a TQFT to
these surfaces. By choosing a TQFT that respects the 4TU relations
(like the standard $X^2=0$ Khovanov TQFT) one can show that it
separates these surfaces (ie, sends them to independent module
maps). The details will not be described here. We do mention that
the reason for $\Psi(4TU)=0$ is coming from the fact that iterative
application of the $3S1$ relation to a 4TU relation is zero as shown
by the picture below (we draw only the pieces of the surfaces that
differ and add the special component on top):

$\eps{10mm}{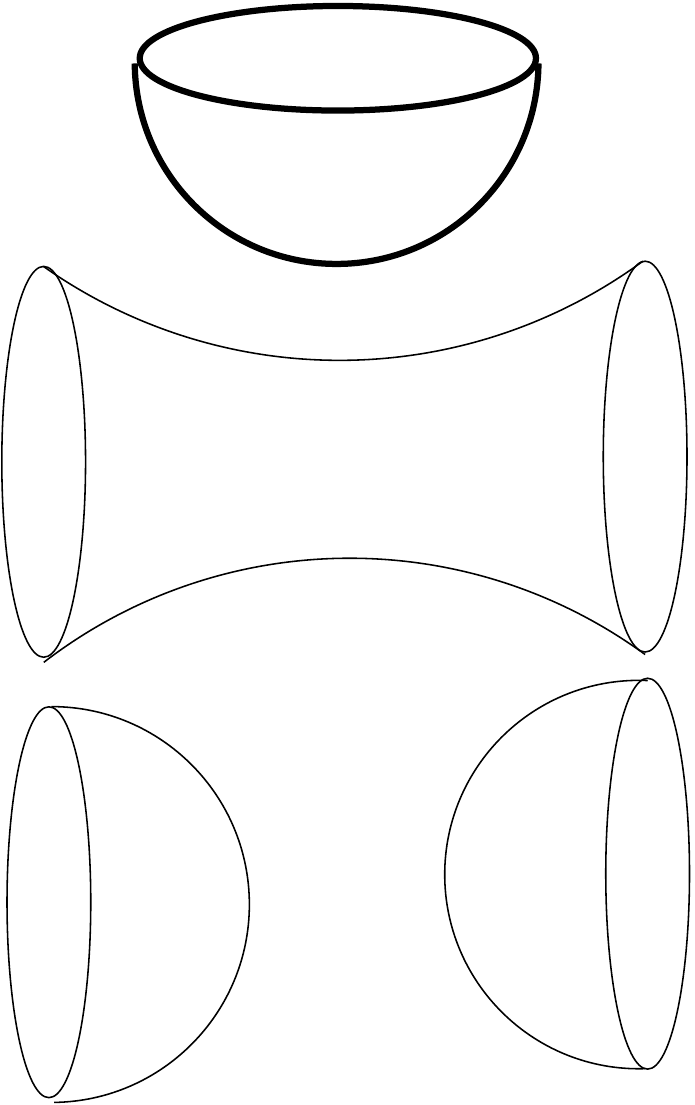} + \eps{10mm}{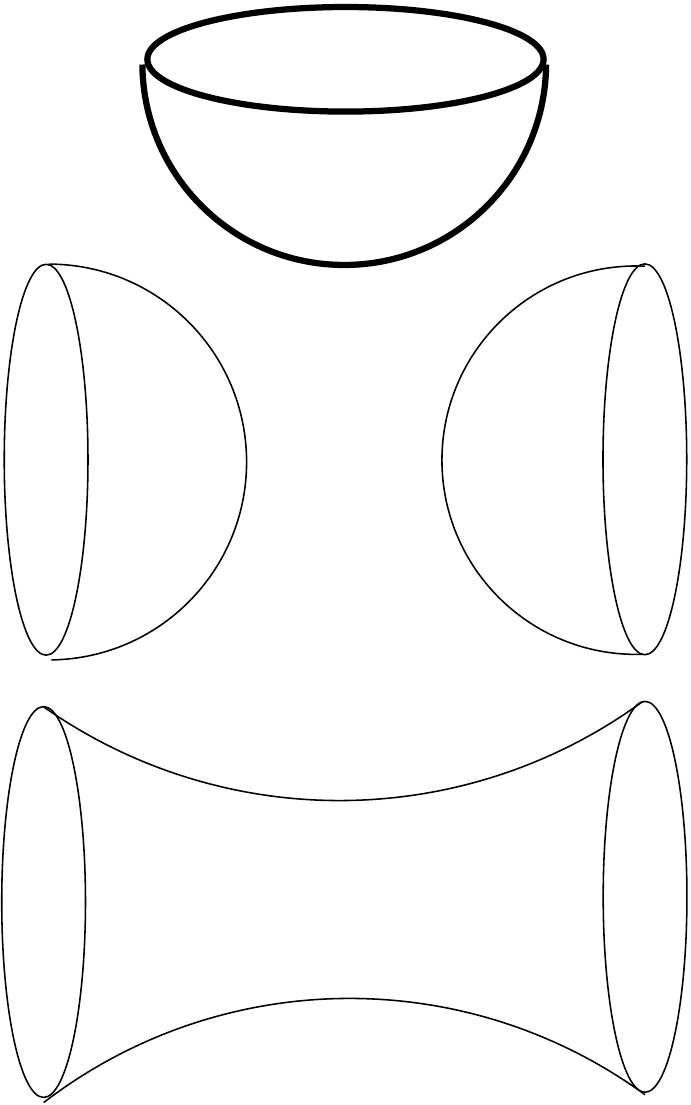} - \eps{14mm}{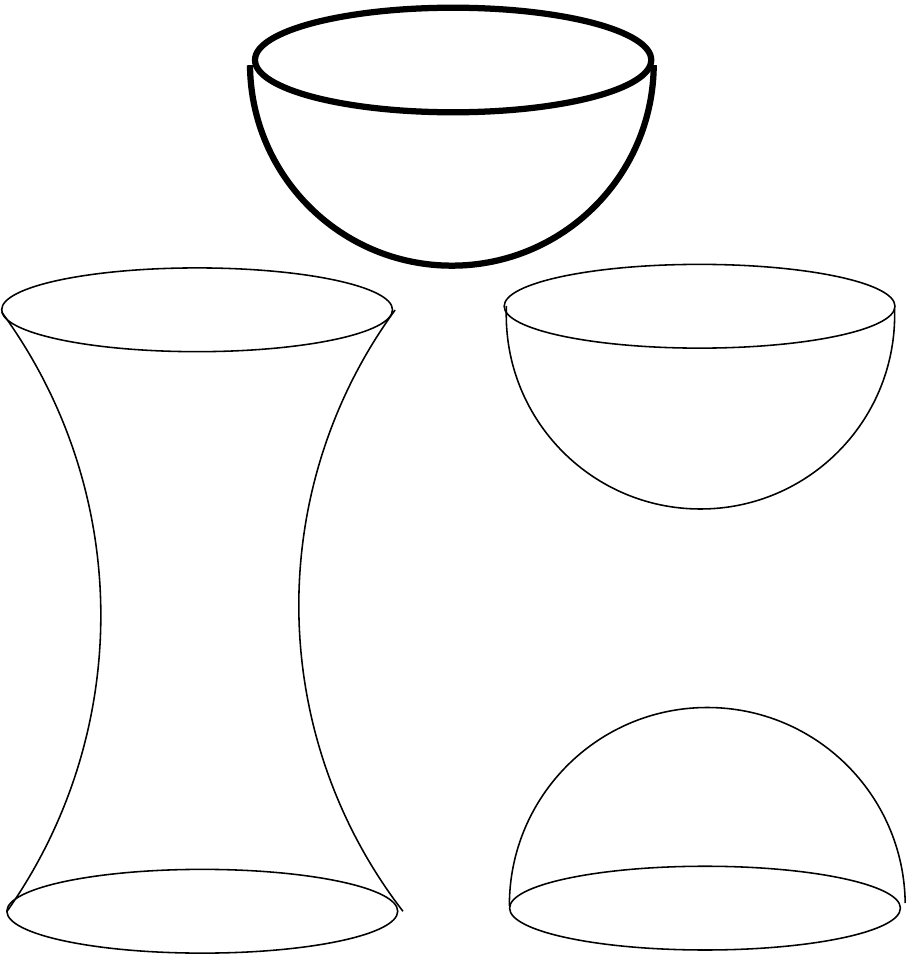}
- \eps{14mm}{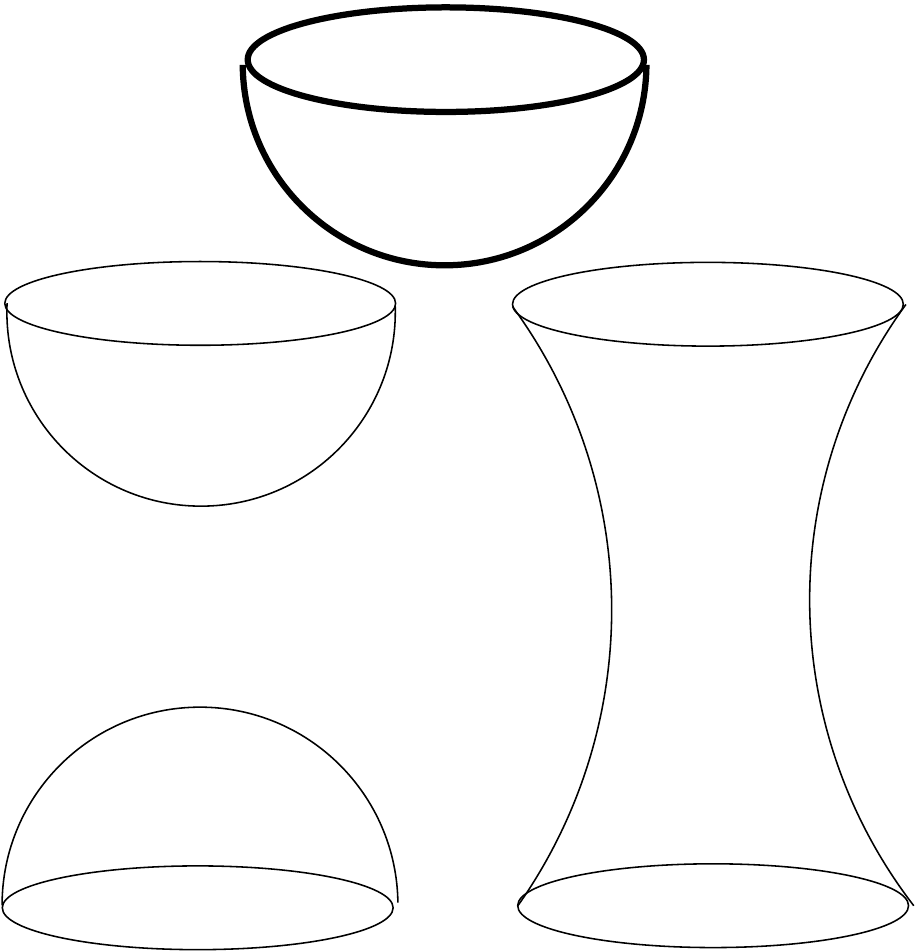} \mapsto \{ \eps{10mm}{} +
\eps{10mm}{} - \eps{10mm}{} \} + \{
\eps{10mm}{} + \eps{10mm}{} - \eps{10mm}{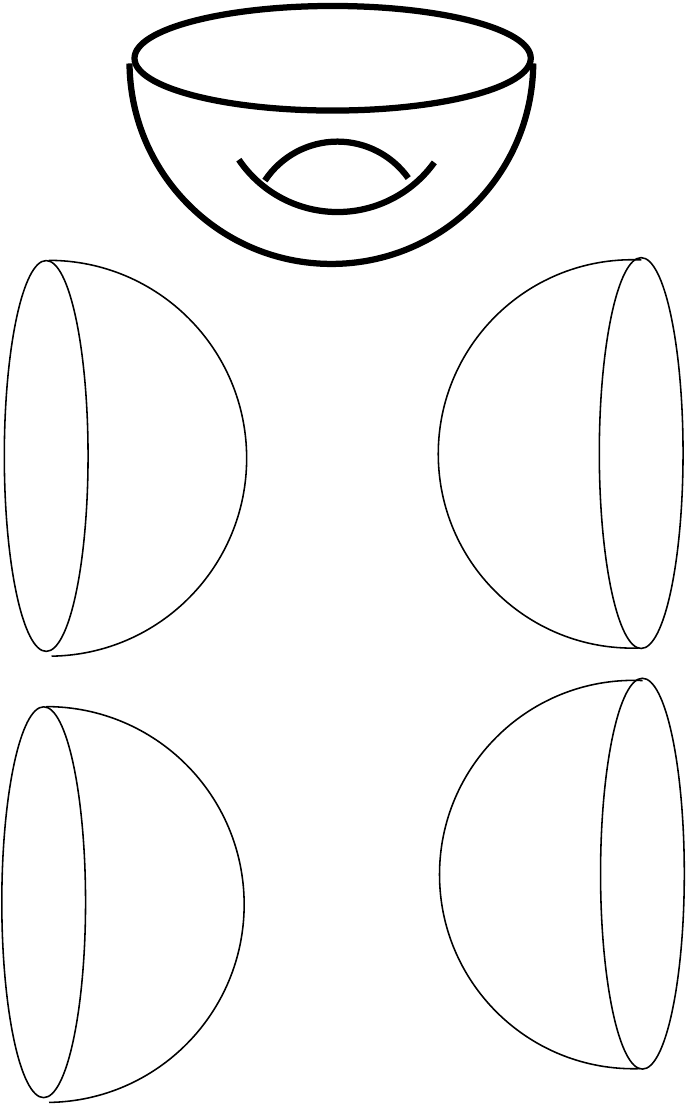}
\} - \{ \eps{10mm}{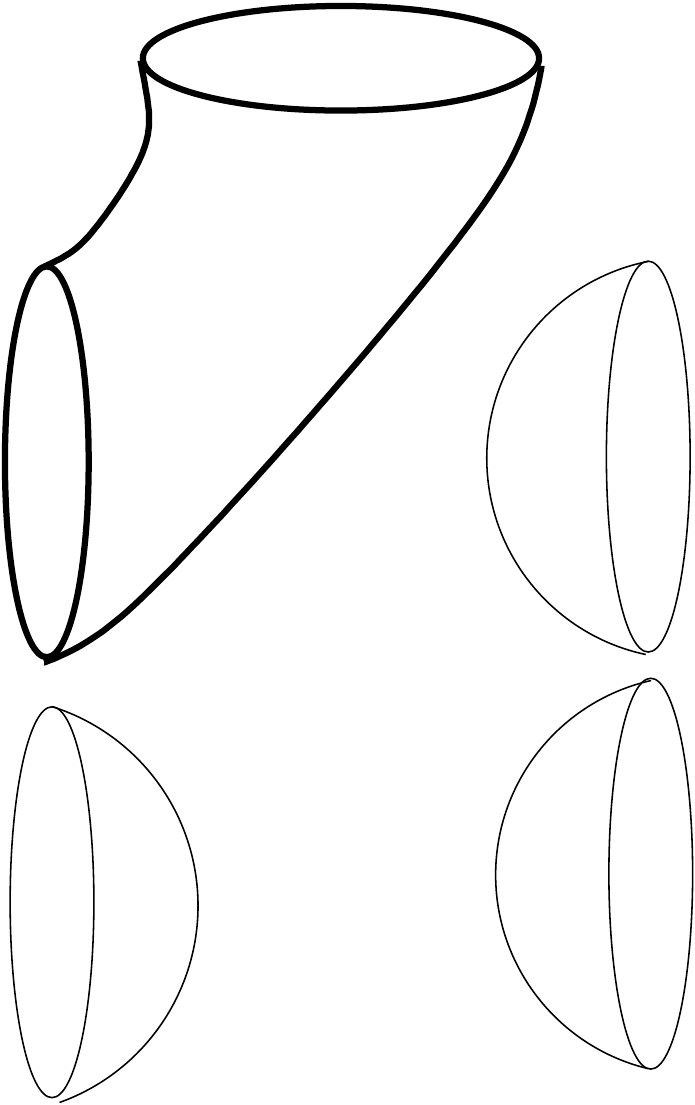} + \eps{10mm}{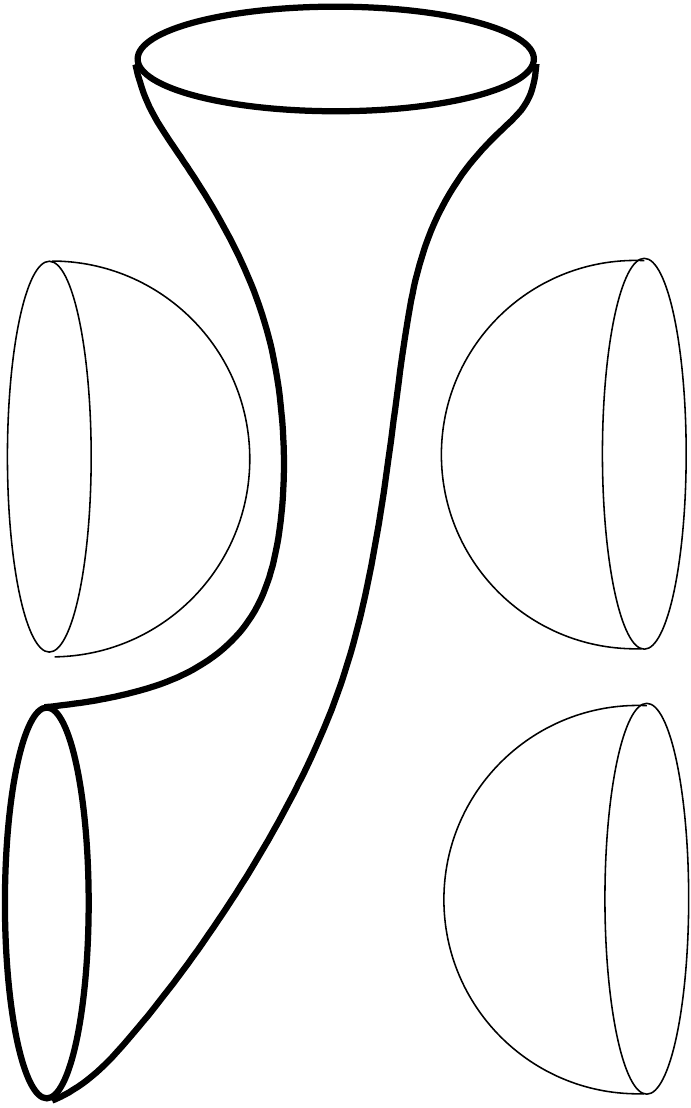} -
\eps{10mm}{welldef7} \} - \{ \eps{10mm}{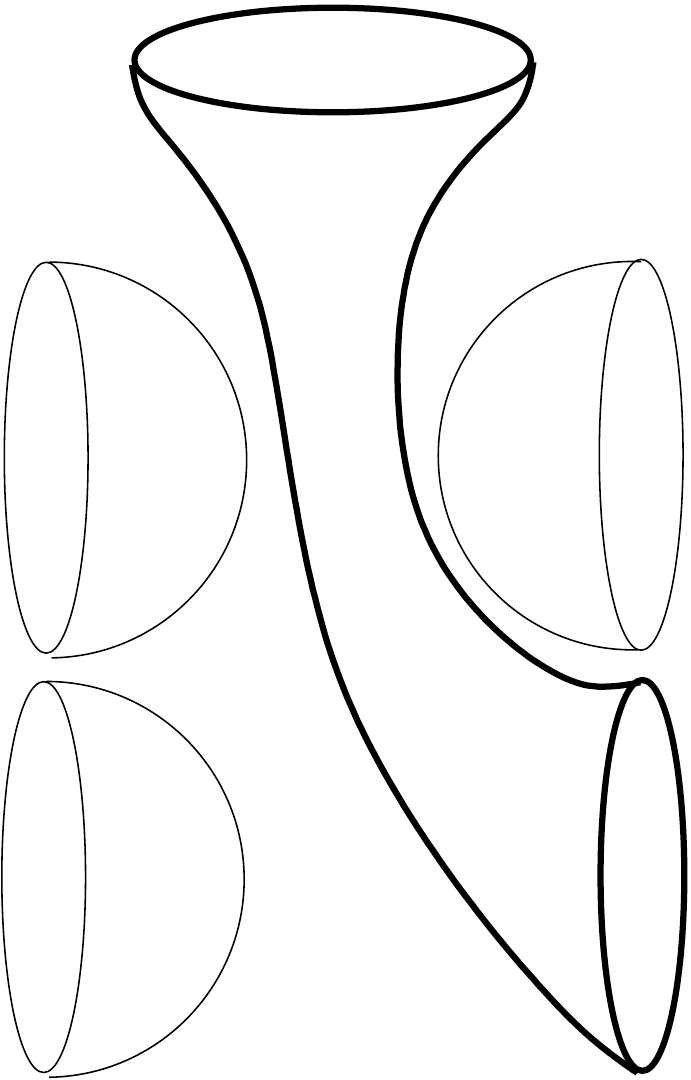} +
\eps{10mm}{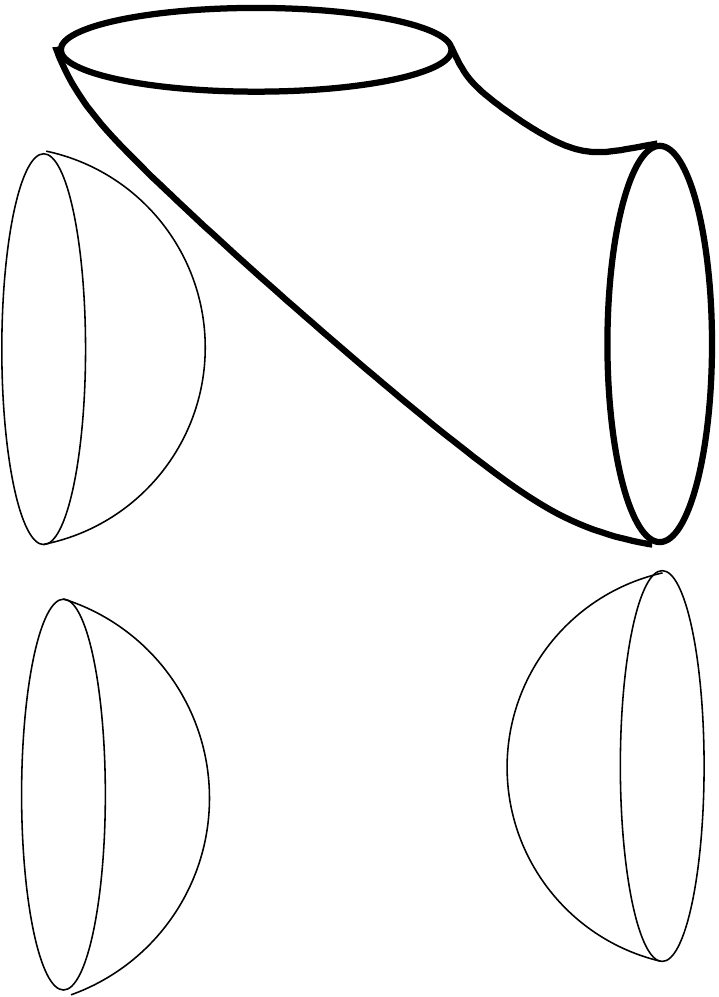} - \eps{10mm}{welldef7} \} = 0 $
\end{proof}

\section{Reduction of the geometric complex associated to a link and
the underlying algebraic structure over $\mathbb{Z}$} 
\label{sec5}

We would like to give a similar treatment for the general case over
$\mathbb{Z}$ as we did over $\mathbb{Q}$. Recall that whenever working
with surfaces over $\mathbb{Z}$ we need to pick a special boundary
circle first. We mark the link at one point (anywhere). After doing
so, we have a special circle at every appearance of an object of
$\Cobl$ in the complex (the one containing the mark). For the sake of
clarity we will always draw the special circle as \emph{a line} (and
call it the special line). This is actually a canonical choice when
working with knots or 1--1 tangles. We defer further discussion on
these issues to \fullref{sec7}. We will use the dotted surface notation,
where a dotted component means that a neck is connected between the
dot and the special component.

\subsection{Reduction theorem over $\mathbb{Z}$}
\vspace{-3pt}

\begin{theorem}\label{Thm2}
The complex associated to a link, over $\mathbb{Z}$, is equivalent
to a complex built from the category of free
$\mathbb{Z}[H]$--modules.
\end{theorem}
\vspace{-3pt}

The theorem is a corollary of the following proposition, which we
will prove first:
\vspace{-3pt}

\begin{proposition}\label{Prop51}
In $\Mat(\Cobl)$ a non-special circle is isomorphic to a column of
two empty set objects. This isomorphism extends (taking degrees into
account) to an isomorphism of complexes in $\Komh\Mat(\Cobl)$ where
the non-special circle objects are replaced by empty sets together
with the appropriate induced complex maps.
\end{proposition}
\vspace{-3pt}

\begin{definition}
The isomorphism from the proposition is called \emph{delooping} and
is given by the following diagram:

\begin{center}
$\eps{18mm}{}$ $\eps{12mm}{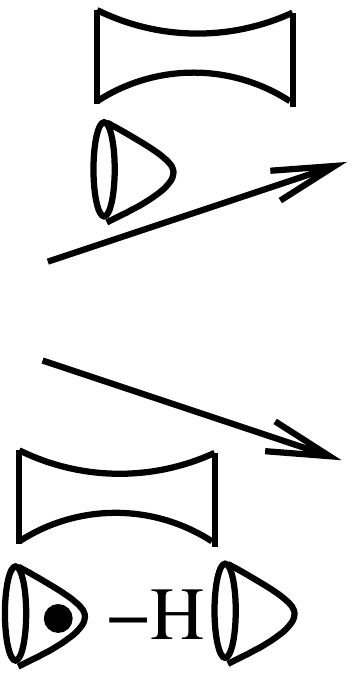}$ $ \left[
\begin{array}{c}
 |~\emptyset \{ -1 \}\\
\\
|~\emptyset\{+1\}
\end{array} \right] $
$\eps{12mm}{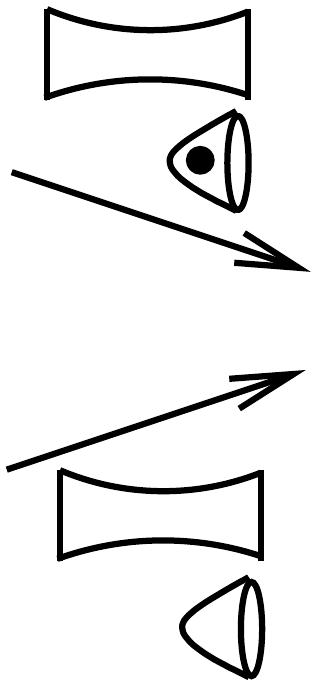}$ $\eps{18mm}{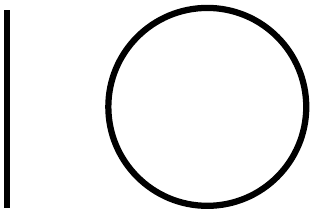}$
\end{center}

This diagram shows how to make a non-special circle disappear into a
column of two empty sets in the presence of the special line. The
special line (remember it is a notation for the special circle)
functions as a ``probe'' which in his presence all the other circles
can be isomorphed into empty sets, leaving us with nothing but the
special line. Note that what used to be a tube between two special
circles is now drawn as \emph{a curtain} between two special lines,
thus a dot will connect a neck to the curtain, and the special
1--handle operator $H$ will add a handle to it.

\end{definition}

\begin{proof}[Proof of \fullref{Prop51}]
The proof is done by looking at the above diagram of morphisms. To
prove that this is indeed an isomorphism we need to trace the arrows
and use the relations of $\Cobl$. When composing the above diagram
from $|~\bigcirc$ to $|~\bigcirc$ one gets: $\eps{10mm}{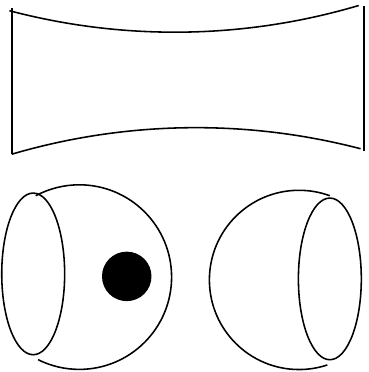} -
\eps{10mm}{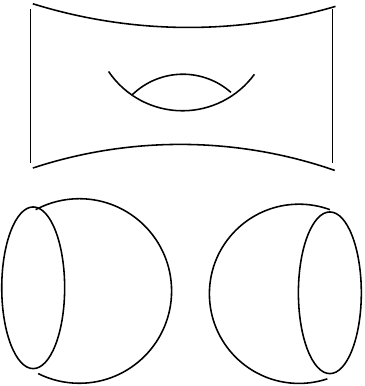} + \eps{10mm}{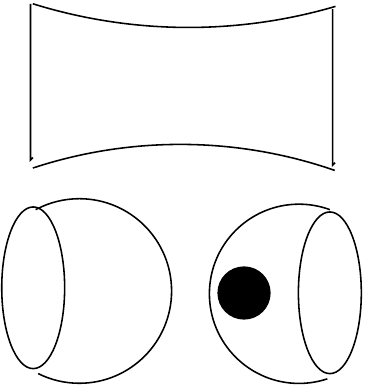}$ (remember that $H$ adds a
handle to the special component). Modulo the $3S1$ relation this
equals to $\eps{10mm}{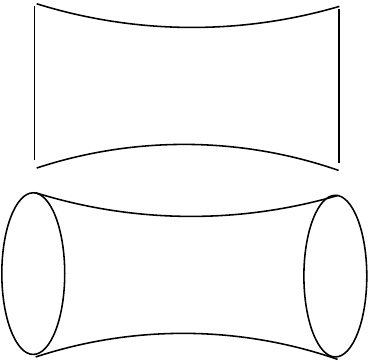}$ which is the identity cobordism. The
other direction is easier, using trivially the S relation
(exercise). The second part of the proposition is a natural
extension of the above diagram to complexes, keeping track of
degrees. Replace the appearances of non-special circles by columns
of empty sets and compose the above isomorphism with the original
complex maps to get the induced maps.
\end{proof}

\begin{proof}[Proof of \fullref{Thm2}]
Given a link, one can mark it at one point and have a special circle
(the one containing the marked point) in each appearance of an
object from $\Cobl$ in the complex associated to the link. Using the
proposition we can reduce this complex and replace all the
non-special circles with columns of empty sets. The only thing left
in the complex are columns with the special line (special circle) in
its entries. The morphism set of that object are curtains (tubes)
with any genus which is isomorphic to $\mathbb{Z}[H]$, where $H$ is
the special 1--handle operator. Thus we get a complex made of columns
of free $\mathbb{Z}[H]$--modules and maps which are matrices with
$\mathbb{Z}[H]$ entries. One would like to check that the functor
from $\Komh\Mat(\Cobl)$ to the category of $\mathbb{Z}[H]$--modules
is indeed well defined (on the morphism level), but this is
tautological from our definitions, as in the proof of the theorem
over $\mathbb{Q}$ (it can still be shown directly by defining the
functor only on generators of $\Cobl$ using our factorization from
the previous section).
\end{proof}

\begin{remark} Though surfaces with two boundary components in
$\Cobl$ have two free generators over $\mathbb{Z}[H]$ only the
connected generator (a curtain) appears in the complex associated to
a link, therefore the morphism group above can be reduced from
$\mathbb{Z}[H] \bigoplus \mathbb{Z}[H]$ to $\mathbb{Z}[H]$. As in
the case over $\mathbb{Q}$ the appearance of $H$ comes only in
homogeneous form, ie, monomials entries, due to grading
considerations.
\end{remark}

\subsection{The algebraic structure underlying the complex
over $\mathbb{Z}$} We follow the same trail as the reduction over
$\mathbb{Q}$ to get some information on the underlying structure of
the geometric complex over $\mathbb{Z}$.

The non-special circle, a basic object in the theory, decomposes
into two copies of another fundamental object (the empty set) with
relative degree shift 2. This restricts functors from the geometric
category to any other category which might carry an algebraic
structure of direct sums (TQFT for example, into the category of
$\mathbb{Z}$--modules). The special line (special circle) stays as
is, but as we will see below it actually carries an intrinsic ``one
dimensional'' object. Later we will also see that the special line
can be \emph{promoted} to carry higher dimensional algebraic
objects.

Denote $|~\{-1\}$ by $v_-$ (comparing to~\cite{ba1}) or X
(comparing to~\cite{kho1}) and $|~\{+1\}$ by $v_+$ or 1 (comparing
accordingly). We use the tensor product symbol to denote unions of
such empty sets (always with one special line only), thus keeping
track of degree shifts.

Now, we take the pair of pants map $\eps{7mm}{InvertedPOP}$ between
$|~\bigcirc\bigcirc$ and $|~\bigcirc$. By this we mean the complexes
$|~\bigcirc\bigcirc\xrightarrow{\rm pants}|~\bigcirc$ and
$|~\bigcirc\xrightarrow{\rm pants}|~\bigcirc\bigcirc$ where the special
lines are connected with a curtain. We reduce these complexes into
empty sets complexes (replacing the non-special circles with empty
sets columns and replacing the maps with the induced maps, as in the
reduction theorem). We denote by $\Delta_1$ the composition:

\begin{center}
$ \Delta_1: \left[ \begin{array}{c}
v_-\\
v_+
\end{array} \right]
\xrightarrow{\rm isomorphism} |~\bigcirc
\xrightarrow{\eps{7mm}{InvertedPOP}} |~\bigcirc\bigcirc
\xrightarrow{\rm isomorphism} \left[
\begin{array}{c}
v_-\otimes v_-\\
v_-\otimes v_+\\
v_+ \otimes v_-\\
v_+ \otimes v_+
\end{array} \right]
$
\end{center}

We denote by $m_1$ the composition in the reverse direction. By
composing surfaces and using the 4TU/S/T relations the maps we get
are (again, these are not linear maps, yet, but a matrix of
cobordisms):
\[
  \Delta_1: \begin{cases}
    v_+ \mapsto  \left[ \begin{array}{c} v_-\otimes v_+ \\ v_+\otimes v_- \\ - H v_+\otimes v_+  \end{array} \right] &
    \\
    v_- \mapsto v_-\otimes v_- &
  \end{cases}
  \qquad
  m_1: \begin{cases}
    v_+\otimes v_-\mapsto v_- &
    v_+\otimes v_+\mapsto v_+ \\
    v_-\otimes v_+\mapsto v_- &
    v_-\otimes v_-\mapsto Hv_-.
  \end{cases}
\]
Let us check what are the induced maps when the pair of pants
$\eps{7mm}{}$ involve the special line. Denote by $\Phi$ the
following composition:

\begin{center}
$ \Phi: \left[ \begin{array}{c}
v_-\\
v_+
\end{array} \right]
\xrightarrow{\rm isomorphism} |~\bigcirc
\xrightarrow{\eps{7mm}{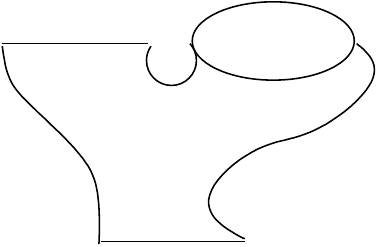}} | $
\end{center}

Denote by $\Psi$ the composition in the reverse direction. Composing
and reducing surfaces we get the following maps:
\[
  \Phi: \begin{cases}
    v_- \mapsto H \\
    v_+ \mapsto 1
  \end{cases}
  \qquad
  \Psi: \begin{cases}
    \emptyset \mapsto v_-
  \end{cases}
\]
($\Delta_1,m_1$) is the type of algebraic structure one
encounters in a TQFT (co-product and product of the Frobenius
algebra). The appearance of a special line (special circle) is most
natural in the construction of the \emph{reduced knot homology},
introduced in~\cite{kho3}. In the setting of Khovanov homology
theory, first, one views the entire chain complex as a complex of
$\mathbb{A}$--modules ($\mathbb{A}$ is the Frobenius algebra
underlying the TQFT used) through a natural action of $\mathbb{A}$
(here one has to mark the knot and encounter the special line).
Then, one can take the kernel complex of multiplication by X to be
the reduced complex. One can also take the \emph{co-reduced
complex}, which is the image complex of multiplication by X. The
above ($\Delta_1$,$m_1$) structure is exactly the TQFT denoted
$\calF_7$ in~\cite{kho1}. The $\mod2$ specialization of this theory
appeared in \cite[Section 9]{ba1}. The additional structure
coming from $\Psi$ and $\Phi$ is exactly the structure of the
\emph{co-reduced} theory of the ($\Delta_1,m_1$) TQFT, where the
special line carries the ``one dimensional'' object generated by X.
It is important to note again that this is an \emph{intrinsic}
information on the underlying structure, coming before any TQFT is
even applied. This gives the above ``pre-TQFT co-reduced structure''
a unique place in the universal theory.
\vspace{-6pt}

\section{TQFTs and link homology theories put on the complex}
\label{sec6}
\vspace{-6pt}

The complex invariant is geometric in nature and one cannot form
homology groups directly (kernels make no sense since the category
is additive but not abelian). Thus we need to apply a functor into
an algebraic category where one can form homology groups. We will
classify such functors and present the universal link homology
theory. We will also explore the relative strength and the
information held within these functors. One of these types of
functors is a TQFT which maps the category $Cob$ into the category
of modules over some ring. In the two dimensional case, which is the
relevant case in link homology, such TQFT structures are equivalent
to Frobenius systems and are classified by them~\cite{ab,kho1}. A
priori, the link homology theory coming from a TQFT does not have to
satisfy the 4TU/S/T relations of the complex.  \fullref{sec63} shows why
any TQFT (up to twisting and base change) used to create link
homology can be put on the geometric complex, making the complex
universal for link homology coming from TQFTs. We show that the
universal link homology functor is actually a TQFT and every other
link homology functor factorizes through it (including non TQFTs).
\vspace{-6pt}

\subsection{Two important TQFTs}
Given the algebraic structures $(\Delta_1, m_1)$ and
$(\Delta_2,m_2)$, introduced in the previous sections, one can add
linearity and get the following algebraic structures (using the same
notation):
\[
  \Delta_1: \begin{cases}
    v_+ \mapsto v_+\otimes v_- + v_-\otimes v_+  - H v_+\otimes v_+ &\\
    v_- \mapsto v_-\otimes v_- &
  \end{cases}
  \]
  \[
  m_1: \begin{cases}
    v_+\otimes v_-\mapsto v_- &
    v_+\otimes v_+\mapsto v_+ \\
    v_-\otimes v_+\mapsto v_- &
    v_-\otimes v_-\mapsto Hv_-
  \end{cases}
\]
\[
  \Delta_2: \begin{cases}
    v_+ \mapsto v_+\otimes v_- + v_-\otimes v_+ &\\
    v_- \mapsto v_-\otimes v_- + T v_+\otimes v_+ &
  \end{cases}
  m_2: \begin{cases}
    v_+\otimes v_-\mapsto v_- &
    v_+\otimes v_+\mapsto v_+ \\
    v_-\otimes v_+\mapsto v_- &
    v_-\otimes v_-\mapsto Tv_+
  \end{cases}
\]
Then, one can construct TQFTs (Frobenius systems) based on these
algebraic structures. The first will be called $\calF_H$ and is
given by the algebra $A_H=R_H[X]/(X^2-HX)$ over $R_H=R[H]$ (system
$\calF_7$ in~\cite{kho1}). The second system is
$A_T=R_T[X]/(X^2-T)$ over the ring $R_T=R[T]$ which will be called
$\calF_T$ (named $\calF_3$ in~\cite{kho1}). Of course one can take
the reduced or co-reduced homology structures, and we denote it by
superscripts ($\calF^{co}_H$ for example). Due to the fact that
these are the underlying algebraic structures of the geometric
complex it is not surprising that they dominate the information in
link homology theory.
\vspace{-6pt}

\subsection{The universal link homology theories}
\vspace{-6pt}

We give the most general link homology theory that can be applied to
the geometric complex.
\vspace{-6pt}

\begin{theorem}
Every functor used to create link homology theory from the geometric
complex over $\mathbb{Q}$, or any ring with 2 invertible, factors
through $\calF_T$. Ie, the TQFT $\calF_T$ is the universal link
homology theory when 2 is invertible and holds the maximum amount of
information.
\end{theorem}

\begin{proof}
Since we showed that in the geometric complex the morphism groups
are $\mathbb{Q}[T]$--modules the target objects of the functor must
be too (by functoriality). Since every object (circle) in the
complex is isomorphic to a direct sum of two empty sets, the functor
will be determined by choosing a single $\mathbb{Q}[T]$--module
corresponding to the empty set. The universal choice would be
$\mathbb{Q}[T]$ itself, and the complex coming from any other module
will be obtained by tensoring the complex with it.
\end{proof}

\begin{theorem} Every functor used to create link homology
theory from the geometric complex over $\mathbb{Z}$ factors through
$\calF^{co}_H$ which holds the maximum amount of information. Ie,
the TQFT $\calF^{co}_H$ is the universal functor for link homology.
\end{theorem}

\begin{proof}
The complex is equivalent to a free $\mathbb{Z}[H]$--modules complex,
and the universal choice is $\mathbb{Z}[H]$. All other theories are
just tensor products of the universal one.
\end{proof}

\subsection{A ``big'' Universal TQFT by Khovanov}\label{sec63}

A priori, TQFTs can be used to create link homology theory without
the geometric complex (meaning without satisfying the 4TU/S/T
relations). In~\cite{kho1} Khovanov presents a universal rank 2
TQFT (Frobenius system) given by the following formula (denoted
$\calF_5$ there):
\[
  \Delta_{ht}: \begin{cases}
    1 \mapsto 1\otimes X + X\otimes 1 -h1\otimes 1&\\
    X \mapsto X\otimes X + t1\otimes 1 &
  \end{cases}
\]
\[
  m_{ht}: \begin{cases}
    1\otimes X\mapsto X &
    1\otimes 1\mapsto 1 \\
    X\otimes 1\mapsto X &
    X\otimes X\mapsto t1 + hX.
  \end{cases}
\]
This structure gives the Frobenius algebra
$A_{ht}=R_{ht}[X]/(X^2-hX-t)$ over the ring
$R_{ht}=\mathbb{Z}[h,t]$. Comparing to our notation (and these of
~\cite{ba1}) is done by putting $v_-=X$ and $v_+=1$. We will denote
this TQFT $\calF_{ht}$.

Given a (tensorial) functor $\calF$ from $Cob$ to the category of
$R$--modules (a TQFT), one can construct a link homology theory and
ask whether it is (homotopy) invariant under Reidemeister moves. If
it is invariant under the first Reidemeister move, then it is a
Frobenius system of rank 2. Khovanov showed that every rank 2
Frobenius system can be \emph{twisted} into a descended theory,
$\calF'$, which is a base change of $\calF_{ht}$ (base change is
just a unital ground ring homomorphism which induces a change in the
algebra. The notion of twisting is explained in~\cite{kho1} and ref
therein). The complexes associated to a link using $\calF$ and
$\calF'$ are isomorphic, thus all the information is still there
after the twist. Since base change just tensors the chain complex
with the appropriate new ring (over the old ring), system
$\calF_{ht}$ is universal for link homology theories coming from
TQFTs. The fact that $\calF_{ht}$ satisfies the 4TU/S/T relations
allows one to apply it as a homology theory functor on the geometric
complex and use~\cite{ba1} results (giving Proposition 6 in
\cite{kho1}). Base change does not change the fact that a theory
satisfies the 4TU/S/T relations, and since $\calF_{ht}$ satisfies
these relations, we have that every tensorial functor $\calF$ that
is invariant under Reidemeister--1 move can be twisted into a functor
$\calF'$ that satisfies the relations of $\Cobl$ and thus can be put
on the geometric complex without losing any homological information
relative to $\calF$. This gives the universality of the geometric
complex for all link homology theories coming from TQFTs. The
universal link homology functor presented in the theorems above
captures all the information coming from TQFT constructions.

Over $\mathbb{Q}$ it is easy to see how our universal theory
captures all the information of khovanov's ``big'' TQFT. We start
with $\calF_{ht}$ over $\mathbb{Q}$ (ie, $R_{ht}=\mathbb{Q}[h,t]$).
By doing a change of basis: $
\begin{cases}
\begin{array}{l}
v_- -\frac{h}{2}v_+ \mapsto v_-\\
v_+ \mapsto v_+
\end{array}
\end{cases}$
we get the theory given by:
\[
  \begin{cases}
    v_+ \mapsto v_+\otimes v_- + v_-\otimes v_+ &\\
    v_- \mapsto v_-\otimes v_- + (t+\frac{h^2}{4}) v_+\otimes v_+ &
  \end{cases}
  \begin{cases}
    v_+\otimes v_-\mapsto v_- &
    v_+\otimes v_+\mapsto v_+ \\
    v_-\otimes v_+\mapsto v_- &
    v_-\otimes v_-\mapsto (t+\frac{h^2}{4})v_+
  \end{cases}
\]
Re-normalizing $t+\frac{h^2}{4}=\frac{T}{4}:=\tilde{T}$, we see that
the above theory is just $\calF_{\tilde{T}}$ with the ground ring
extended by another superficial variable $h$. Every Calculation done
using $\calF_{ht}$ to get link homology, is equal to the same
calculation done with $\calF_{\tilde{T}}$ tensored with
$\mathbb{Q}[h]$, and holds exactly the same information. Note that
our complex reduction is performing this change of basis (to ``kick
out'' one redundant variable) on the complex level, before even
applying any TQFT. Doing first the complex reduction and then
applying the TQFT $\calF_{ht}$ will factorize the result through
$\calF_{\tilde{T}}$.

In the general case, over $\mathbb{Z}$, one needs some further
observations which we turn to now.

\subsection{About the 1--handle and 2--handle operators in
$\calF_{ht}$} Given the TQFT $\calF_{ht}$ we can ask what is the
1--handle operator $H$ of this theory. Meaning, we want to know what
is the operator that adds a handle to a cobordism, or put in other
words, what is $\eps{5mm}{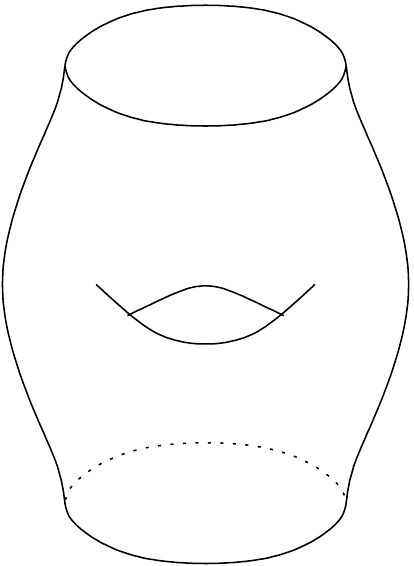}$ viewed as a map $A \mapsto A$
between the two Frobenius algebras associated to the boundary
circles. One can easily give the answer by looking at the
multiplication and co-multiplication formulas: $H=2X-h$. For the
readers who know a bit about topological Landau--Ginzburg models, $H$
is the Hessian of the potential, and indeed the results match . The
2--handle operator $T$ follows immediately by computing $T=H^2$ and
reducing modulo $X^2-hX-t$ to get: $T=4t+h^2$. The 2--handle operator
is an element of the ground ring, as expected from our 2--handle
lemma, and can be multiplied anywhere in a tensor product. The
1--handle operator is an element of $A$ but not of the ground ring,
and thus when operating on a tensor product one needs to specify the
component to operate on, as we did by picking the special circle.
Since $A$ is two dimensional, picking the basis to be $\left(
\begin{smallmatrix} 1 \\ X \end{smallmatrix} \right)$ would give a two
dimensional representation of the operators above: $H=\left(
\begin{smallmatrix} -h&2t \\ 2&h \end{smallmatrix} \right)$, $T=\left(
\begin{smallmatrix} 4t+h^2 & 0 \\ 0&4t+h^2  \end{smallmatrix} \right)$.

\subsection{The special line and H promotion}
The main result of the previous section was a reduction of the
complex associated to a link into a complex composed of columns of
the special line and maps which are matrices with $H$ monomials
entries. Given such a complex one wants to create a homology theory
out of it, ie, apply some functor to it that will put an $H$--module
on the special line with the possibility of taking kernels (and
forming homology groups). As we have seen, the intrinsic structure
that the special line caries is the one dimensional module
$\mathbb{Z}[H]$ generated by $X$. We can get a different structure
by replacing $H$ with any integer matrix of dimension $n$ and the
special line with direct sum of $n$ copies of $\mathbb{Z}$. We call
this type of process \emph{a promotion} (realizing it is a fancy
name for tensoring process). Another type of promotion is replacing
$H$ with a matrix (of dimension $n$) with polynomial entries in two
variables ($h$ and $t$, say) and then promoting the special line
into the direct sum of $n$ copies of $\mathbb{Z}[h,t]$. Note that
some promotions lose information.

\subsection{TQFTs via promotion and new type of link homology functors}

First, we want to find a promotion of $H$ that is equivalent to the
theory $\calF_{ht}$ and which does not lose any information on the
complex level. From the above it is clear
that such a promotion is $H=\left( \begin{array}{cc} -h&2t \\
2&h \end{array} \right)$. The special line will be promoted into two
copies of $\mathbb{Z}[h,t]$. Since each power of the promoted $H$
adds a power of $t$ and $h$ to the entries, the power of the
1--handle operator $H$ can be uniquely determined after the promotion
and we lose no information. The algebraic complex one gets by
applying this promotion is equal to the complex one gets by applying
the TQFT $\calF_{ht}$ before the complex reduction.

We can now easily get the other familiar TQFTs using the same
complex reduction and $H$ promotion technique:
\[
\calF_{\tilde{T}} \leftarrow H=\left( \begin{array}{cc} 0&2\tilde{T} \\
2&0 \end{array} \right) \hspace{7mm} \calF_H \leftarrow H=\left( \begin{array}{cc} -H&0 \\
2&H \end{array} \right)
\]
One gets the standard Khovanov homology ($X^2=0$) by the promotion
$H=\left( \begin{array}{cc} 0&0 \\ 2&0 \end{array} \right)$. The
special line is promoted to double copies of
$\mathbb{Z}[\tilde{T}]$, $\mathbb{Z}[H]$ and $\mathbb{Z}$
respectively. One can get the \emph{reduced Khovanov homology} by
focusing on the top left entry of the $H$ promotion matrix. It is
interesting to note that in the standard Khovanov homology case
indeed we lose information and $H^2=T=0$. In the geometric language
this means that in order to get the standard Khovanov homology one
has to ignore surfaces with genus 2 and above. This was observed in
\cite{ba1}. Reduction to Lee's theory~\cite{lee1} is done by
substituting $\tilde{T}=1$.

One can create other types of promotions that will enable us to
control the order of $H$ involved in the theory. We are able to
cascade down from the most general theory, the one that involves all
powers of $H$ (ie, surfaces with any genera in the topological
language), into a theory that involves only certain powers of $H$
(ie, genera up to a certain number). For example, promote the
special line to three copies of $\mathbb{Z}$, and $H$ to the matrix
$\left(
\begin{array}{ccc} 0&0&0 \\ 1&0&0 \\ 0&1&0 \end{array} \right)$. This
theory involves only powers of $H$ smaller or equal to 2, that is
surfaces of genus up to 2. These promotions can be viewed as a
family of theories extrapolating between the ``standard'' Khovanov
TQFT and our universal theory (reminding of a perturbation expansion
in string theory).

\section{Further discussion}\label{sec7}

\subsection{Comments on the question of homotopy classes
versus homology\\theories} The geometric complex associated to a link is
an invariant up to homotopy of complexes. Thus its fullest strength
lies in the homotopy class of the complex itself, and have the
potential of being a stronger invariant than any functor applied to
the complex to produce a homology theory. The reduction given above
might be the beginning of an approach to the following question:
\emph{classify all complexes associated to links up to homotopy}.
The category $\Komh\Mat(\Cobl)$ seemed at first too big and
complicated for an answer, but complexes built on free modules over
polynomial rings look more hopeful. The isomorphism of complexes
reduce that question into the following question:

{\sl
\medskip 
{\bf Question 1}\qua Classify all homotopy types of chain complexes
with free $\mathbb{Z}[H]$--modules as chain groups and $\mathbb{Z}[H]$
matrices as maps (without loss of generality, with monomial entries).}

\medskip 
Another interesting question regarding the strength of
the complex invariant is the following: 

\textsl{Do homology theories
(functors from the topological category to an algebraic one)
completely classify homotopy classes of complexes?} 

Combined with question 1 we get:

{\sl
\medskip 
{\bf Question 2}\qua Does the geometric complex contain more
information than all the possible algebraic complexes (homology
theories) that can be put on it? In other words, are there two
non-homotopic geometric complexes, associated to links in the
geometric formalism, that are not separated by some functor to an
algebraic category?}

{\sl
\medskip 
{\bf Question 3}\qua  If the answer is no, then do all the possible
homology groups of these algebraic complexes classify the homotopy
class of the geometric complex?.}

The answers will determine the relative strength of the geometric
complex invariant. The complex reduction answers question 2.

\medskip
{\bf Answer 2}\qua {\sl No}

{\bf Explanation}\qua The complex is built from a
category equivalent to a category of free  $\mathbb{Z}[H]$--modules.
Thus, it has a faithful algebraic representation, ie, we have found
a homology theory that represents faithfully the complexes in
$\Komh\Mat(\Cobl)$. Moreover, this homology theory is the co-reduced
theory of the TQFT ($\Delta_1,m_1$), which can be reached by
applying a specific tautological functor. The geometric complex
holds the same amount of information coming from the chain complex
of this specific homology theory (TQFT).

We are left then with:

{\sl
\medskip 
{\bf Question 3 (revised)}\qua  Do the homology groups of complexes in
$\Komh\Mat(\mathbb{Z}[H])$, associated to links, classify the
homotopy type of the chain complexes?}

\subsection{Comments on marking one of the boundary circles}
\vspace{-6pt}

We classified all the surfaces with at least one boundary circle
over $\mathbb{Z}$. These are the ones relevant for link homology. In
order to do that we picked up a special circle and marked it. The
presentation of the generators depends on which circle we choose,
but this choice has no importance for the classification itself and
the topology of the generating surfaces. If the link is a knot, then
in the context of knot homology and the geometric complex there is a
natural way of determining a special circle in each and every
appearance of an object of $\Cobl$ in the geometric complex. This is
done by marking a point on the knot we start with --- the special
circle will be the circle with the mark on it. Marking the knot is
not a new procedure in knot invariants theory and appears in many
parts of quantum invariants theory. As far as this part of knot
theory is concerned the marking of the knot has no effect on theory.
One can look at this process as marking a point on the knot for
cutting it open to a 1--1 tangle (or a ``long knot'') --- the theories
of knots and long knots (1--1 tangles) are ``isomorphic'' in our
context. Once we have a 1--1 tangle, there is always a special line
appearing naturally in the complex. When we deal with links, one
might choose different components of the link to place the mark but
the choice does not matter and gives isomorphic complexes (they
might be presented differently though). This is obvious from the
fact that the complex reduction is local, and thus in every
appearance of an object of $\Cobl$ in the complex one can choose the
special circle \emph{independently} and apply the complex
isomorphisms. Different choices are linked through a series of
complex isomorphisms.
\vspace{-6pt}

\subsection{Comments on computations}
\vspace{-6pt}
As was shown in
\cite{ba2}, for the case of the original standard Khovanov
homology, fast computations can make one happy! At first sight the
complex isomorphism presented above does not seem to reduce the
geometric complex at all (it doubles the amount of objects in the
complex). The surprising thing is that one can use this isomorphism
(and some homology reduction technics) to create a very efficient
crossing-by-crossing local algorithm to calculate the geometric
complex and link homology. This was done for the standard Khovanov
homology (high genera are set to zero) as reported in
\cite{ba2,ba4}, and more recently done for the universal case over
$\mathbb{Z}[H]$, as reported in~\cite{ba3} (J Green implementing
the algorithm by D Bar-Natan based on the work of the author
presented here). It is important to mention that this is currently
the fastest program to calculate Khovanov homology (and the complex
itself), and that the isomorphism of the geometric complex (the
delooping process) is a crucial component of it~\cite{ba4}.
Computations of the universal complex (and related issues, like
various promotions of that complex) will be treated in a future
work.
\vspace{-10pt}

\subsection{Comments on tautological functors}
\vspace{-10pt}
In \cite[Section 9]{ba1} tautological functors where defined on
the geometric complex associated to a link. One needs to fix an
object in $\Cobl$, say $\calO'$, and then the tautological functor
is defined by $\calF_{\calO'}(\calO)=\Mor(\calO',\calO)$, taking
morphisms to compositions of morphisms. Our classification allows us
to state the following regarding tautological functors over the
geometric complex:
\vspace{-10pt}

\begin{corollary}
The tautological functor $\calF_\bigcirc(-)$ over $\mathbb{Z}[H]$ is
the TQFT $\calF_H$.

The tautological functor $\calF_\emptyset(-)$ over $\mathbb{Q}[T]$
is the TQFT $\calF_T$.
\end{corollary}
\vspace{-10pt}

\begin{proof}
Indeed this is a corollary to the surface classification. In the
first case declare the special circle to be the source $\bigcirc$.
\end{proof}
\vspace{-10pt}

Since every TQFT (after twisting) is factorized through
$\calF_H/\calF_T$, tautological functors hold all the information
one can get. Moreover, since the information held in the geometric
complex is manifested in the theories $\calF_H/\calF_T$, it seems
that asking about homology theories which are not tautological is
not important. Also it seems, that there is no need in asking about
non-tensorial functors (ie, functors for which
$\calF(\bigcirc\bigcirc)\ncong\calF(\bigcirc)\otimes\calF(\bigcirc)$).
For example, the tautological functor $\calF_{\bigcirc^n}(-)$ is
equivalent to $\calF_\bigcirc(-) \otimes
\calF_\bigcirc(\bigcirc)^{\otimes n-1}$ and thus holds the same
information as the the ones in the corollary, which can be
considered as universal for that sake. The question whether every
functor on $\Cobl$ can be represented as a tautological functor
seems also less important due to the above.
\vspace{-10pt}

\subsection{Comments on embedded versus abstract surfaces} As
noticed in \cite[Section 11]{ba1}, when one looks at surfaces in
$\Cobl$ over a ground ring with the number 2 invertible there is no
difference between working with embedded surfaces (inside a cylinder
say) or with abstract surfaces. This is true due to the fact that
any knotting of the surface can be undone by cutting necks and
pulling tubes to unknot the surface. In other words, by cutting and
gluing back, using the NC relation (divided by 2) both ways, one can
go from any knotted surface to the unknotted version of it embedded
in 3 dimensional space. Our claim is that the same is true even when
2 is not invertible. The proof is a similar argument applied to any
knotted surface using the $3S1$ relation:

\begin{center}
$\begin{array}{ccccccc} \eps{20mm}{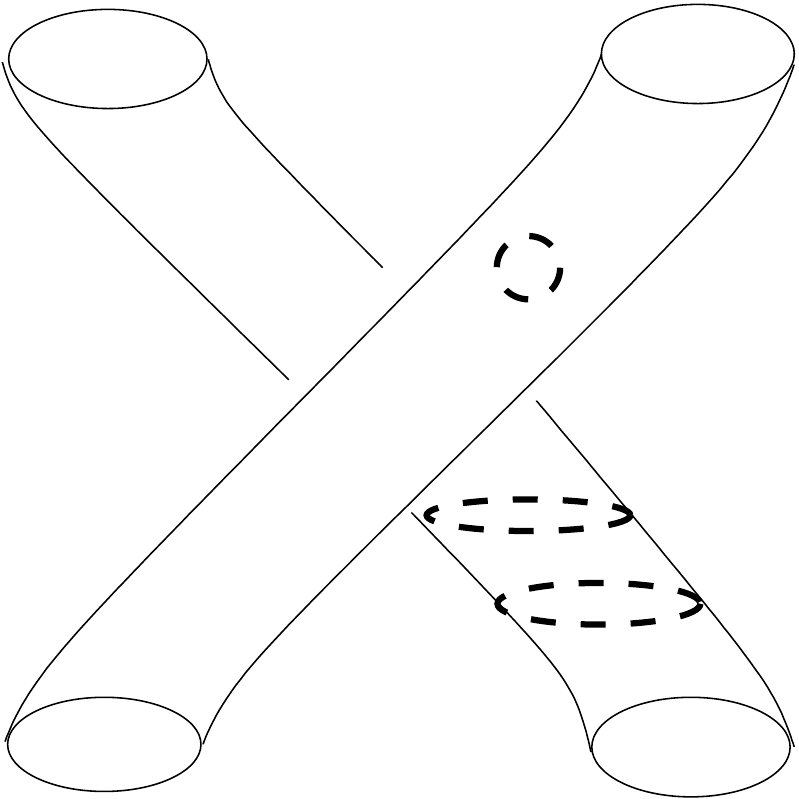} & = & \eps{20mm}{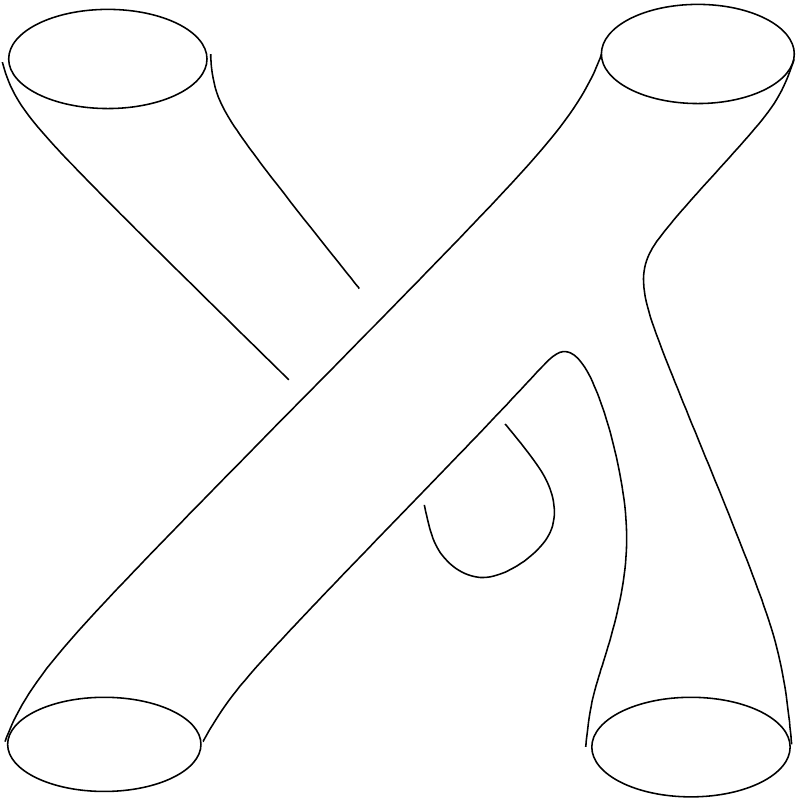} & +
& \eps{20mm}{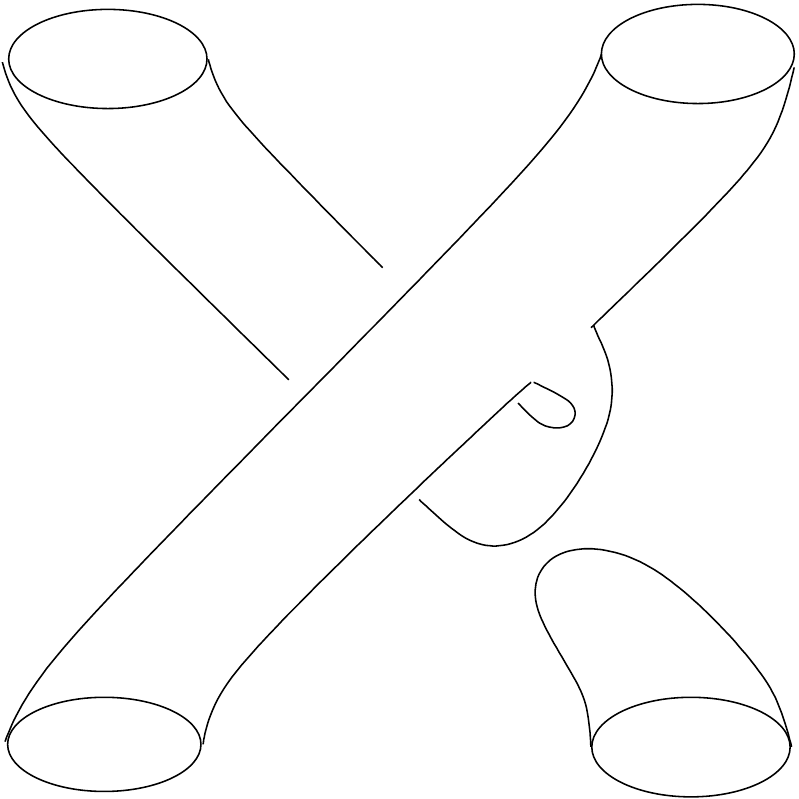} & - & \eps{20mm}{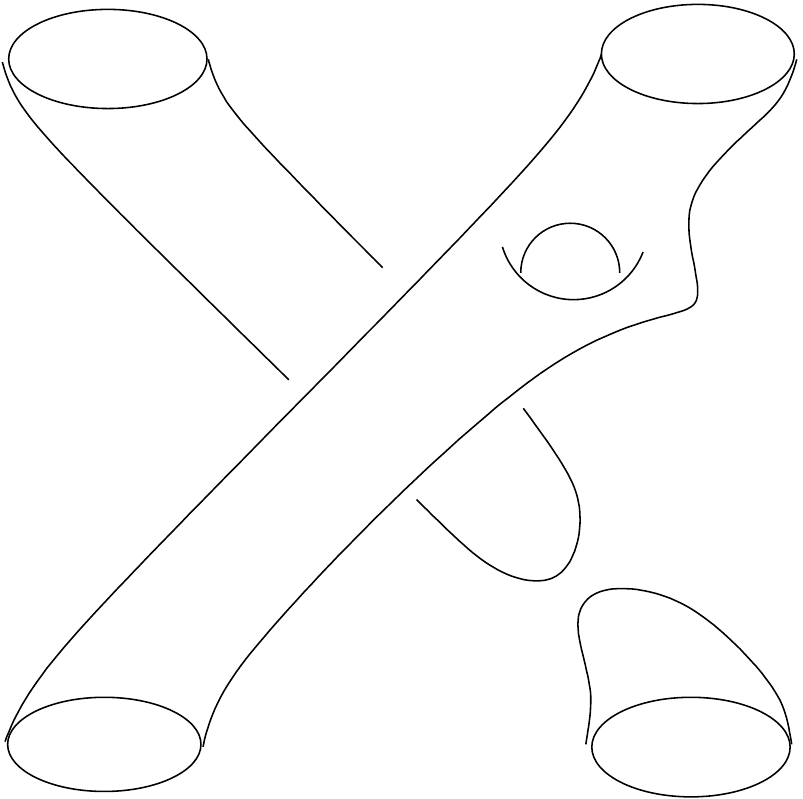} \\
 & & \updownarrow & & \updownarrow & & \updownarrow\\
 \eps{20mm}{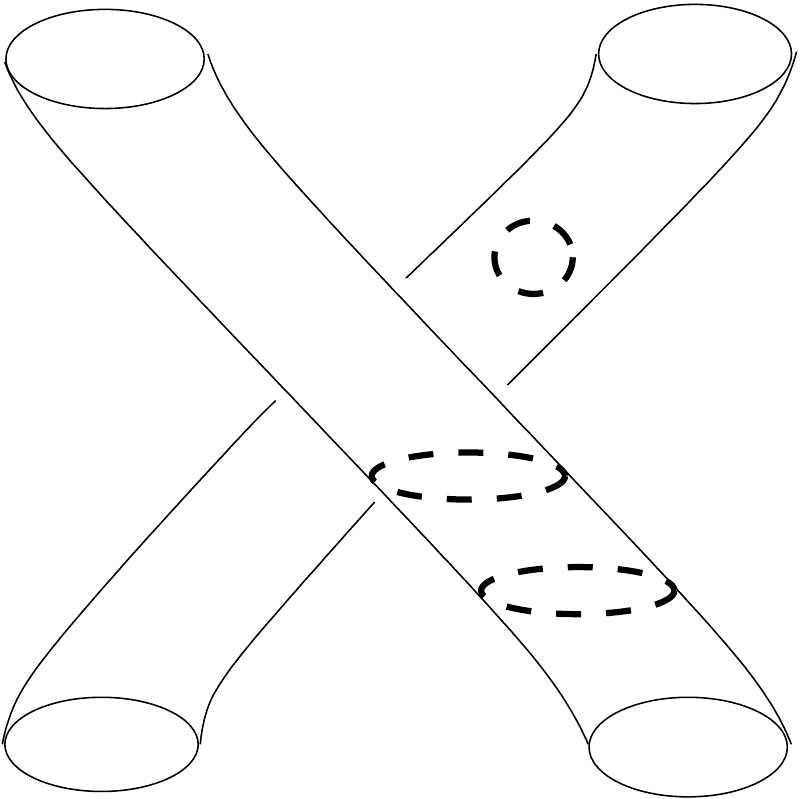} & = & \eps{20mm}{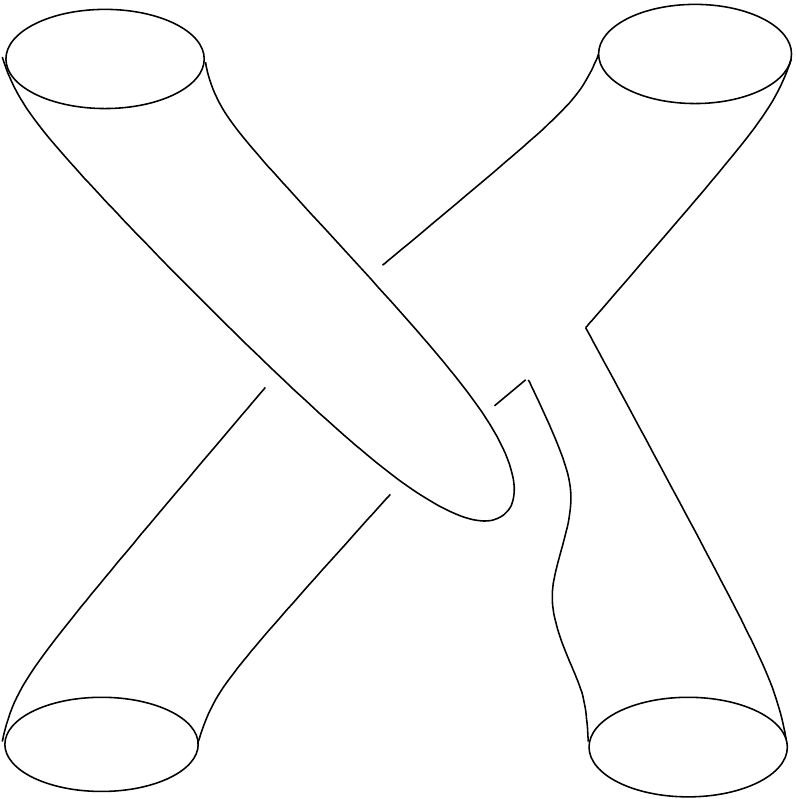} & + & \eps{20mm}{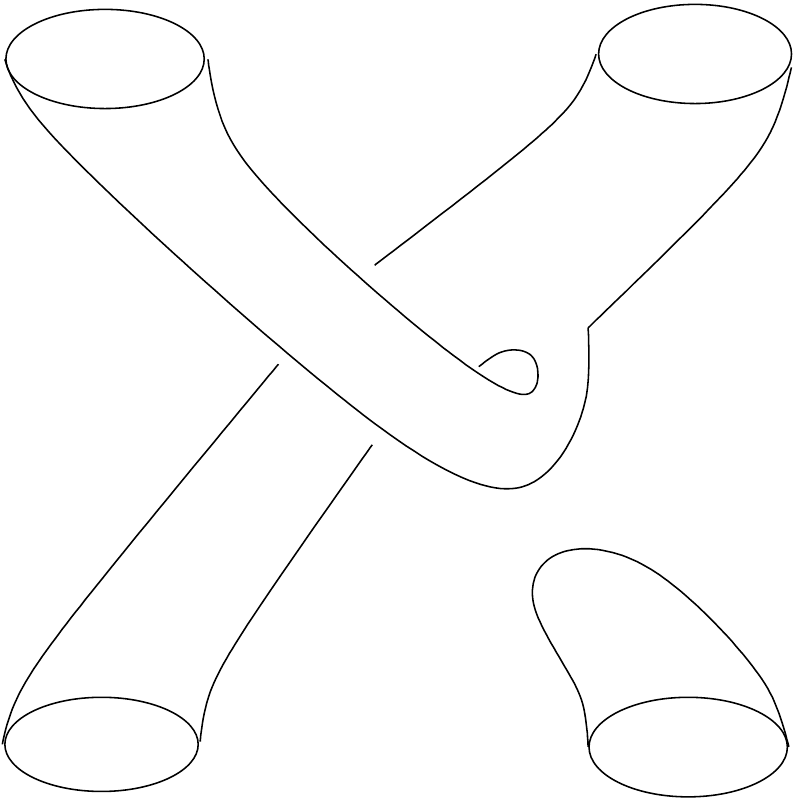} & 
- & \eps{20mm}{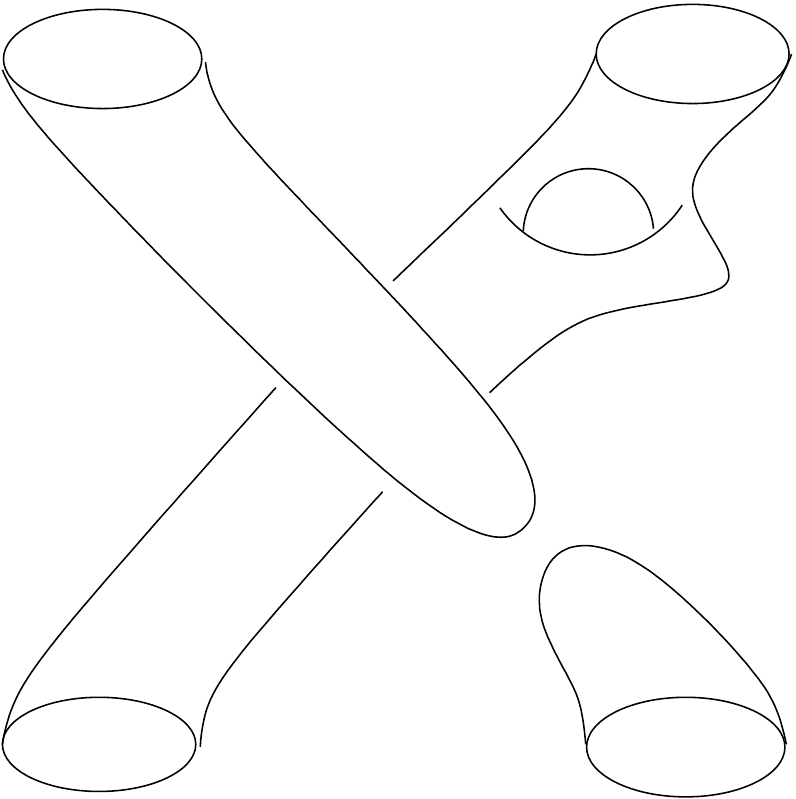}
 \end{array}$
\end{center}

The above picture shows that every crossing (a part of a knotted
surface embedded in 3 dimensions) can be flipped using the $3S1$
relation twice. Apply the $3S1$ relation once on the dashed sites
(going from top left), then smoothly change the surface (going down
the arrows) and finally use the dashed sites for another application
of the $3S1$ relation (reaching the bottom left). Every embedded
surface can be unknotted this way, justifying previous comments
about ignoring the issues of embedding (\fullref{sec2}).

\bibliographystyle{gtart}
\bibliography{link}

\begin{thebibliography}{}
\providecommand\bibmarginpar{\leavevmode\marginpar}
\def\urlstyle#1{{\tt #1}}

\bibitem{ab}
\textbf{L Abrams}, \href{http://dx.doi.org/10.1142/S0218216596000333}
  {\emph{Two-dimensional topological quantum field theories and {F}robenius
  algebras}}, J. Knot Theory Ramifications 5 (1996) 569--587 \xox{MR}{1414088}

\bibitem{ba3}
\textbf{D Bar-Natan}, \emph{4TU-2}, talk given at UQAM, Montreal (October 2005)
\ Available at \setbox0\hbox{\makeatletter\@url
{http://www.math.toronto.edu/~drorbn/Talks/}}
\href{http://www.math.toronto.edu/~drorbn/Talks/}
{\unhbox0}

\bibitem{ba4}
\textbf{D Bar-Natan}, \emph{Fast Khovanov Homology Computations}
  \xox{arXiv}{math.GT/0606318}

\bibitem{ba2}
\textbf{D Bar-Natan}, \emph{I've Computed Kh(T(9,5)) and I'm Happy}, talk given
  at George Washington University (February 2005)
\ Available at \setbox0\hbox{\makeatletter\@url
{http://www.math.toronto.edu/~drorbn/Talks/}}
\href{http://www.math.toronto.edu/~drorbn/Talks/}
{\unhbox0}

\bibitem{ba1}
\textbf{D Bar-Natan}, \href{http://dx.doi.org/10.2140/gt.2005.9.1443}
  {\emph{Khovanov's homology for tangles and cobordisms}}, Geom. Topol. 9
  (2005) 1443--1499 \xox{MR}{2174270}

\bibitem{kho2}
\textbf{M Khovanov},
  \href{http://projecteuclid.org/getRecord?id=euclid.dmj/1092749199} {\emph{A
  categorification of the {J}ones polynomial}}, Duke Math. J. 101 (2000)
  359--426 \xox{MR}{1740682}

\bibitem{kho3}
\textbf{M Khovanov},
  \href{http://projecteuclid.org/getRecord?id=euclid.em/1087329238}
  {\emph{Patterns in knot cohomology I}}, Experiment. Math. 12 (2003) 365--374
  \xox{MR}{2034399}

\bibitem{kho1}
\textbf{M Khovanov}, \emph{Link homology and {F}robenius extensions}, Fund.
  Math. 190 (2006) 179--190 \xox{MR}{2232858}

\bibitem{lee1}
\textbf{E\,S Lee}, \emph{On Khovanov invariant for alternating links}
  \xox{arXiv}{math.GT/0210213}

\end{thebibliography}

\end{document}